\documentclass[12pt,a4paper]{article}
\usepackage[latin1]{inputenc}
\usepackage[english]{babel}
\usepackage{amsmath}                    
\usepackage{amssymb}                    
\usepackage{amsthm}                     
\usepackage{multicol}                   
\usepackage{enumerate}		
\usepackage[pdftex]{graphicx}
\usepackage{epstopdf}			
\usepackage{epsfig} 
  \let\8=\infty \let\0=\emptyset

   \def\Im{\mathop{\rm Im }\nolimits}

  \def\epsilon{\varepsilon}

  \def\N{\mathbb{N}} \def\Z{\mathbb{Z}} 
  \def\R{\mathbb{R}} \def\C{\mathbb{C}}

  \def\F{\mathcal F} 
  \def\G{\mathcal G} 

  \newtheorem{defi}{Definition}[section]
  \newtheorem{teo}{Theorem}
  \newtheorem{pro}[defi]{Proposition}
  \newtheorem{cor}[defi]{Corollary}

  \newtheorem{conj}{Conjecture}
  \newtheorem{OP}{Open Problem}

 \newsavebox\boxdem
 \newif\ifdem \demtrue
 \newenvironment{dem}%
  {\setbox\boxdem\vbox\bgroup\par\medskip\noindent{\underline{Proof:}}}%
  {\hfill$\square$\par\bigskip\egroup\ifdem\unvbox\boxdem\else\fi}
\title{On forbidden sets}
\author{Francisco Balibrea and Antonio Cascales}
\date{\vspace{-5ex}}

\begin{document}
\maketitle\footnote{This paper will be published in the special ECIT'14 issue of Journal of Difference Equations and Applications}
\begin{center}
{\scriptsize
Universidad de Murcia, Departamento de Matemáticas\\
Campus de Espinardo, 30100, Murcia (Spain)\\
balibrea@um.es, antoniocascales@yahoo.es}
\end{center}
\begin{abstract}
In recent literature there are an increasing number of papers where the forbidden sets of difference equations are computed. We review and complete different attempts to describe the forbidden set and propose new perspectives for further research and a list of open problems in this field.
\end{abstract}
\section{Introduction}
The study of difference equations (DEs) is an interesting and useful branch of discrete dynamical systems due to their variety of behaviors and ability to modelling phenomena of applied sciences (see \cite{dynamics,elaydi,periodicities,dynamics2order} and the references therein). The standard framework for this study is to consider iteration functions and sets of initial conditions in such a way that the values of the iterates belong to the domain of definition of the iteration function and therefore the solutions are always well defined. For example, in rational difference equations (RDEs) a common hypothesis is to consider positive coefficients and initial conditions also positive, see \cite{Camouzis2009,dynamics,dynamics2order}.\\
Such kind of restrictions are also motivated by the use of DE as applied models, where negative initial conditions and/or parameters are usually meaningless, \cite{sedaghat2003nonlinear}.\\

But there is a recent interest to extend the known results to a new framework where initial conditions can be taken to be arbitrary real numbers and no restrictions are imposed to iteration functions. It is in this setting where appears the forbidden set of a DE, the set of initial conditions for which after a finite number of iterates we reach a value outside the domain of definition of the iteration function. Indeed, the central problem of the theory of DEs is reformulated in the following way:\\

\emph{Given a DE, to determine the good ($\G$) and forbidden ($\F$) sets of initial conditions. For points in the good set, to describe the dynamical properties of the solutions generated by them: boundedness, stability, periodicity, asymptotic behavior, etc.}\\

In this paper we are interested in the first part of the former problem: how can the forbidden set of a given DE of order $k$ be determined? In the previous literature to describe such sets, when it is achieved, is usually interpreted as to be able to write a general term of a sequence of hypersurfaces in $\R^k$. But those cases are precisely the corresponding to DE where it is also possible to give a general term for the solutions. Unfortunately there are a little number of DEs with explicitly defined solutions. Hence we claim that new qualitative perspectives must be assumed to deal with the above problem. Therefore, the goals of this paper are the following: to organize several techniques used in the literature for the explicit determination of the forbidden set, revealing their resemblance in some cases and giving some hints about how they can be generalized. Thus we get a long list of DEs with known forbidden set that can be used as a frame to deal with the more ambitious problem of describe the 
forbidden 
set of a general DE. We review and introduce some methods to work also in that case. And finally we propose some future directions of research.\\

The paper is organized as follows: after some preliminaries, we review the Riccati DE, which is one of the few examples of DE where the former explicit description is possible. As far as we know, almost all the literature where the forbidden set is described using a general term includes some kind of semiconjugacy with a Riccati DE. DEs obtained via a change of variables or topological semiconjugacy are the topic of the rest of the section. In the following we will discuss how algebraic invariants can be used to transform a given equation into a Riccati or linear one depending upon a parameter, and therefore determining its forbidden set.

After we will deal with an example, found in \cite{TheForbiddenSet_CDV}, of description where the elements of the forbidden set are given recurrently but in explicit form.

We introduce a symbolic description of complex and real points of $\F$ in section \ref{sec:Symbolic}, whereupon  in section \ref{sec:Qualitative} we study some additional ways to deal with the forbidden set without an explicit formula. We finalize with a list of open problems and conjectures.\\

To avoid an overly exposition, we have omitted some topics as the study of systems of difference equations (SDEs), (see \cite{Bajo2014,2014BajoInvariantQuadrics,BajoFrancoPeran2011,BajoLiz2006,GroveEtAl2001,1999_MagnuckaB,2011StevicOnASystem,2012StevicOnSomeSolvable,YazlikTolluTaskara2013,YazlikTolluTaskara2013_Fibonacci}), the use of normal forms (see \cite{RubioMassegu2007}), the systematic study of forbidden sets in the globally periodic case and the importance of forbidden sets in Lorenz maps.

\section{The forbidden set problem}
In this paper we deal with difference equations (DEs) and systems of difference equations (SDEs). General definitions of these concepts are the following.\\
Let $\mathbb{X}$ be a nonempty set and $A\subset \mathbb{X}$ a nonempty subset of $\mathbb{X}$. Let $f:A\rightarrow \mathbb{X}$ be a map. A DE of order $1$ is the formal expression
\begin{equation}
x_{n+1}=f(x_n)
 \label{eq:DEorder1}
\end{equation}
(\ref{eq:DEorder1}) represents a set of finite or infinite sequences of $A$, as given $x_0\in A$, $(x_n)_{n=0}^{M}$ is constructed by recurrence using (\ref{eq:DEorder1}), that is, $x_1=f(x_0)$ and if $x_{n}\in A$ then $x_{n+1}=f(x_{n})$. When this process can be repeated indefinitely then $M=+\8$ and $(x_n)_{n=0}^{M}$ is named the solution of (\ref{eq:DEorder1}) generated by the initial condition $x_0$.\\
We call $\mathcal{F}=\{x_0\in A: \exists n\geq0 \:|\: f(x_n)\notin A\text{ or } f(x_n) \text{ is not defined}\}$ the \emph{forbidden set (FS)} of (\ref{eq:DEorder1}), and $\mathcal G=A\setminus\mathcal F$ \emph{the good set (GS)} of (\ref{eq:DEorder1}).\\
A DE of order $k$ is
\begin{equation}
x_{n+1}=f(x_n,\ldots,x_{n-k+1})
 \label{eq:DEorderk}
\end{equation}
where now $f:A\subset \mathbb{X}^k\rightarrow \mathbb{X}$ is the iteration map, and the FS is defined as
$$\begin{array}{lll}
\mathcal F & =&\{(x_0,\ldots,x_{-k+1}): \exists n\geq0 \:|\: (f(x_n,\ldots,x_{n-k+1}),x_n,\ldots,x_{n-k})\notin A\\
&&\text{ or }f(x_n,\ldots,x_{n-k+1})\text{ is not defined}\}
\end{array}$$
In a similar way, we define a SDE of order $1$ as
\begin{equation}
\left.
\begin{array}{ccc}
  x^1_{n+1}	&=&	f_1(x^1_n,\ldots,x^r_n)\\
  \vdots		&&	\vdots\\
  x^r_{n+1}	&=&	f_r(x^1_n,\ldots,x^r_n)
\end{array}
\right\}
 \label{eq:SDEorder1}
\end{equation}
provided $f_i:A\subset \mathbb{X}^r\rightarrow \mathbb{X}$, $i=1,\ldots,r$.\\
System (\ref{eq:SDEorder1}) can be expressed as DE of order $1$ of type (\ref{eq:DEorder1}) using the vectorial notation $X_n=(x^1_n,\ldots,x^r_n)$ and considering the map $F:A\subset \mathbb{X}^r\rightarrow \mathbb{X}^r$ whose components are $f_1,\ldots,f_r$.\\
Finally, a SDE of order $k$ is defined as a set of equations
\begin{equation}
x^i_{n+1}=f_i(x^1_n,\ldots,x^1_{n-k+1},\ldots,x^r_n,\ldots,x^r_{n-k+1})
 \label{eq:SDEorderk}
\end{equation}
using maps $f_i:A\subset \mathbb{X}^{rk}\rightarrow \mathbb{X}$ for $i=1,\ldots,r$.\\
In vectorial form the set of equations (\ref{eq:SDEorderk}) can be rewritten as in (\ref{eq:DEorderk})
$$X_{n+1}=F(X_n,\ldots,X_{n-k+1})$$
with $F:A\subset\mathbb{X}^{rk}\rightarrow \mathbb{X}^r$.\\

The former definitions depend on how is given the domain of definition $A$ of the iteration map. To remark that point, let's consider $c\in\R\setminus\{0\}$  and the real DE
$$x_{n+1}=c$$
which forbidden set is $\mathcal F=A$ when $A$ is defined as $\R\setminus\{c\}$ meanwhile $\mathcal F$ is empty when $A$ is the natural domain of the constant function.\\

To avoid degenerated or trival cases, some further restrictions must be imposed to $A$. Natural and common restrictions consist usually in regarding which is the domain of definition of the iteration map. Therefore, in equation
\begin{equation}
x_{n+1}=\frac1{x_n}
 \label{eq:1overx}
\end{equation}
we say that $\mathcal F=\{0\}$ as every solution with $x_0\neq0$ is 2-periodic and $x_1$ is not defined when $x_0=0$. It is implicitly assumed that $A=\R\setminus\{0\}$.\\
We remark that (\ref{eq:1overx}) is also a DE over $\mathbb X=\mathbb{RP}$, the projective line, and then $\mathcal F=\emptyset$ being every solution 2-periodic if $A=\mathbb X$ using the rules $\frac10=\8$, $\frac1{\8}=0$.\\

Analogously, the DE
\begin{equation}
x_{n+1}=\frac1{x_n^2+1}
 \label{eq:1overx^2+1}
\end{equation}
has $\mathcal F=\emptyset$ in $\R$, $i\in \mathcal F$ when (\ref{eq:1overx^2+1}) is taken over the complex field and $\mathcal F$ is again empty in the projective plane $\mathbb{CP}$.\\

In practical applications an undefined zero division means that the denominator belongs to certain neighbourhood of zero. Therefore in DE as
\begin{equation}
 x_{n+1}=\frac{P(x_n)}{Q(x_n)}
 \label{eq:RDEorder1}
\end{equation}
where $P$, $Q$ are real polynomials, the forbidden set problem could be studied using $A=\R\setminus(-\varepsilon,\varepsilon)$ for a machine-value $\varepsilon>0$.\\

In applied models where only positive values of the variables have practical meaning, $A$ could be $(\R^+)^k$.\\

As a final example of different ways to consider the $A$ set, let's recall the DE associated to a Lorenz map. Let $I\subset\R$ be the unit interval $[0,1]$. Let $c\in (0,1)$. We say that $f:I\setminus\{c\}\rightarrow I$ is an expansive Lorenz map if the following conditions are satisfied
\begin{enumerate}[i)]
 \item $f$ is increasing and continuous in each interval $[0,c)$ and $(c,1]$
 \item $\lim\limits_{x\rightarrow c^-}f(x)=1$ and $\lim\limits_{x\rightarrow c^+}f(x)=0$
 \item $f$ is topologically expansive, i.e., there exists $\varepsilon>0$ such that for any two solutions $(x_n)_{n=0}^{+\8}$ and $(y_n)_{n=0}^{+\8}$ of (\ref{eq:DEorder1}) not containing the point $c$, there is some $i\geq 0$ with $|x_i-y_i|>\varepsilon$, and if some of the former solutions contains $c$ the inequality remains valid taking $f(c)=1$ or $f(c)=0$.
\end{enumerate}
Condition (iii) is equivalent to say that the preimages of the point $c$ are dense in $I$, or, in the forbidden set notation, to say that $\overline{\mathcal F}=I$.\\
It is obvious that $f$ can be arbitrarily defined in $c$ (without bilateral continuity) and in that case $\mathcal F=\emptyset$. In the standard definition it is assumed that $A=I\setminus\{c\}$.\\
Lorenz maps are an important tool in the study of the Lorenz differential equations and the Lorenz attractor, and also in the computation of topological entropy in real discontinuous maps. See \cite{1996GlendinningHall_ZerosOfTheKneading,1990_HubbardSparrow} and the references therein.\\

In the following we deal mostly with DEs where iteration functions are quotients of polynomials, known as rational difference equations (RDEs).
\section{Semiconjugacies on Riccati equations}
Let's briefly recall some well known results about Riccati DEs. Let $a,b,c,d\in\R$ such that $|c|+|d|\neq0$. A Riccati DE of order $1$ is
\begin{equation}
x_{n+1}=\frac{a+bx_n}{c+dx_n}
 \label{eq:RiccatiOrder1}
\end{equation}
Some special cases occur for particular values of the parameters. If $d=0$, $ad-cb=0$ or $b+c=0$ the equation is linear, constant or globally of period $2$ respectively (see \cite{tesisACV,dynamics2order}). If none of those conditions stand, an affine change of variables transforms (\ref{eq:RiccatiOrder1}) into
\begin{equation}
y_{n+1}=1-\frac R{y_n}
 \label{eq:RiccatiOrder1_SF}
\end{equation}
$R=\frac{bc-ad}{(b+c)^2}$ being the Riccati number. Finally $y_n=\frac{z_n}{z_{n-1}}$ leads to the linear DE $z_{n+1}=z_n-Rz_{n-1}$. The closed form solution of the linear DE ables to write a closed form solution for (\ref{eq:RiccatiOrder1_SF}). Moreover, there is a correspondence between the solutions of the linear DE containing the zero element and the finite or forbidden set generated solutions of (\ref{eq:RiccatiOrder1_SF}). This is the idea behind the characterization of the FS of (\ref{eq:RiccatiOrder1}) that can be found, for example, in \cite{AziziThesis,tesisACV,dynamics2order}.

From a topological point of view, that characterization shows that $\mathcal F$ is a convergent sequence or a dense set in $\R$ or a finite set. Moreover, in the last case equation (\ref{eq:RiccatiOrder1_SF}) is globally periodic (see \cite{RubioMassegu2009}).\\

Analogously, for the second order Riccati DE
\begin{equation}
x_{n+1}=a+\frac b{x_n}+\frac c{x_n x_{n-1}}~~~a,b,c\in\R,~~c\neq0
 \label{eq:RiccatiOrder2}
\end{equation}
we can transform it into a linear one using the change of variables $x_n=\frac{y_n}{y_{n-1}}$. The forbidden set is then described as (see \cite{Azizi2012,AziziThesis,RiccatiOrderTwo_Sedaghat})
\begin{equation}
\mathcal F  =\bigcup_{n=-1}^{+\8}\left\{
(u,v)\in\R^2 \:|\: \beta_{1n}uv+\beta_{2n}u+\beta_{3n}=0
\right\}
 \end{equation}
where coefficients $\left\{\beta_{1n},\beta_{2n},\beta_{3n}\right\}_{n=-1}^{+\8}$ depend on the roots of the characteristic equation. Roughly speaking $\mathcal F$ is a countable union of plane hyperbolas convergent to a limit curve in some cases, dense in open sets in other cases and even a finite collection when the DE is globally periodic.\\
Our first proposal (open problem \ref{op:QualitativeDescriptionRiccati}) is to clarify which kind of topological objects can be obtained in the Riccati DE  of second order and to generalize to the Riccati DE of order $k$
\begin{equation}
x_{n+1}=a_0+\frac{a_1}{x_n}+\frac{a_2}{x_nx_{n-1}}+\ldots+\frac{a_k}{x_n\ldots x_{n-k+1}}~~~a_0,\ldots,a_k\in \R,~~a_k\neq0
 \label{eq:RiccatiOrderk}
\end{equation}
We claim also (conjecture \ref{conj:RiccatiGloballyPeriodic}) that DE (\ref{eq:RiccatiOrderk}) is globally periodic if and only if $\mathcal F$ is a finite collection of hypersurfaces in $\R^k$, generalizating what happens in the case of order 1 (see \cite{RubioMassegu2009}).

To study the higher order Riccati DE see also \cite{Azizi2013}.\\

In the former discussion, it was important to determine the set of initial conditions whose solutions in the linear equation included the zero element. This is called the zero set of the linear equation, $\mathcal Z$. For example, in the Riccati DE of order $1$, $\mathcal Z$ is defined as the set of points $(z_0,z_{-1})\in\R^2$ such that the solution of the linear DE includes the zero element. But $z_0$ and $z_{-1}$ depend on the initial condition $y_0$ of (\ref{eq:RiccatiOrder1_SF}) via the change of variables $y_n=\frac{z_n}{z_{n-1}}$. We can write the dependence as $z_0=z_0(y_0)$, $z_{-1}=z_{-1}(y_0)$. Therefore, $\mathcal F$ can be expressed as
$$\begin{array}{ccl}
\mathcal F	&=&\{y_0\in\R:\exists n\geq0 \:|\: z_n(y_0)=0\}   \\
		&=&\{y_0\in\R:(z_{0}(y_0),z_{-1}(y_0))\in\mathcal Z\}
  \end{array}
$$
Generalizing to the Riccati DE of order $k$ (\ref{eq:RiccatiOrderk}), we get
$$\begin{array}{ccl}
\mathcal F	&=&\{(x_{0},\ldots,x_{-k+1})\in\R^k:\exists n\geq0 \:|\: z_n(x_0,\ldots,x_{-k+1})=0\}\\
		&=&\{(x_{0},\ldots,x_{-k+1})\in\R^k:(z_{0},\ldots,z_{-k})\in\mathcal Z\}
  \end{array}
$$
where $z_n=z_n(x_0,\ldots,x_{-k+1})$ are the terms of the associated linear equation of order $k+1$.\\
Note that the semiconjugacy formula $y_n=\frac{z_n}{z_{n-1}}$ has zero as a pole. The forbidden can be regarded as the transformation of certain special set of the linear equation corresponding to the singularities of the semiconjugacy.\\

There is an increasing number of works in the recent literature that use the closed form solution of the Riccati equation to describe the forbidden set of other rational difference equations. Let's review and complete some of them.\\

The order $3$ DE
\begin{equation}
x_{n+1}=\frac{ax_nx_{n-1}}{bx_n+cx_{n-2}}
 \label{eq:AboZeid_1}
\end{equation}
can be transformed into the linear form $z_{n+1}=\frac{b}a+\frac c a z_n$ using the semiconjugacy $z_n=\frac{x_{n-2}}{x_n}$. The forbidden set is a sequence of planes in $\R^3$. The case $b=1$ and $a,c>0$ in (\ref{eq:AboZeid_1}) was developed in \cite{2009_Sedaghat_GlobalBehaviours}, while $a=\pm 1$, $b=1$ and $c=-1$ is in \cite{2014_AboZeid_TwoThirdOrder}. See also \cite{AboZeid2014,2014AboZeid_3order} for $a>0$ and $bc<0$.\\

A similar change of variables, $z_n=\frac{x_{n-3}}{x_n}$, gives the former linear equation when it is applied to the following DE of order $4$, see \cite{2014AboZeid_4order}.
\begin{equation}
x_{n+1}=\frac{ax_nx_{n-2}}{bx_n+cx_{n-3}}~~~a,b,c>0
\label{eq:AboZeid_2}
\end{equation}
We propose a generalization of this problem remarking that the following DE
\begin{equation}
x_{n+1}=\frac{ax_nx_{n-k+1}}{bx_n+cx_{n-k}}~~~a\neq0, |b|+|c|>0
\label{eq:AboZeid_generalized}
\end{equation}
becomes $z_{n+1}=\frac{b}a+\frac c a z_n$ with $z_n=\frac{x_{n-k}}{x_n}$ (open problem \ref{op:AboZeid_generalized}).\\

In \cite{BajoLiz2011}, the change $z_n=x_nx_{n-1}$ gives the Riccati equation $z_{n+1}=\frac{z_n}{a+bz_n}$ when we apply it to
\begin{equation}
x_{n+1}=\frac{x_{n-1}}{a+bx_nx_{n-1}}
 \label{eq:BajoLiz2011}
\end{equation}
The authors explicitly describe $\mathcal F$ as a sequence of plane hyperbolas.\\

Equation (\ref{eq:BajoLiz2011}) admits several generalizations. If we use the change of variables $z_n=x_nx_{n-k}$ applied to the Riccati DE $z_{n+1}=\frac{z_n}{a+bz_n}$, we obtain
\begin{equation}
x_{n+1}=\frac{x_nx_{n-k}}{ax_{n-k+1}+x_nx_{n-k+1}x_{n-k}} 
\label{eq:BajoLiz2011_generalized1}
\end{equation}
a RDE of order $k+1$ whose forbidden set can be described using that of the Riccati DE. In this case we have a RDE of third degree. A more interesting family of equations, where the degree remains equal to two, is
\begin{equation}
 x_{n+1}=\frac{x_{n-2i-1}}{a+bx_{n-i}x_{n-2i-1}}
\label{eq:BajoLiz2011_generalized2}
\end{equation}
where $i$ is a natural number. Remark that for $i=0$ we get equation (\ref{eq:BajoLiz2011}). Using the change of variables $z_n=x_{n-i-1}x_n$, (\ref{eq:BajoLiz2011_generalized2}) leads to
$$z_{n+1}=\frac{z_{n-i}}{a+bz_{n-i}}$$
that it is not a Riccati difference equation, but can still be reduced to a linear form. Indeed, if with an affine change we transform it into
$$y_{n+1}=1-\frac R{y_{n-i}}$$
we can introduce $y_n=\frac{u_n}{u_{n-i-1}}$ to get
$$u_{n+1}=u_{n-i}-Ru_{n-2i-1}$$
a linear DE of order $2i+2$ (see open problem \ref{op:BajoLiz2011_generalized}).\\

In \cite{Mcgrath2006}, it is shown that $\F$ is a family of straight lines in $\R^2$ in the case of equation
\begin{equation}
x_{n+1}=\frac{ax_{n-1}+bx_n}{cx_{n-1}+dx_n}x_n
 \label{eq:Mcgrath2006}
\end{equation}
Here we get the Riccati form $z_{n+1}=\frac{a+bz_n}{c+dz_n}$ via the change of variables $z_n=\frac{x_n}{x_{n-1}}$.\\
Given $i\geq0$ a possible generalization is 
\begin{equation}
x_{n+1}=\frac{ax_{n-2i-1}+bx_{n-i}}{cx_{n-2i-1}+dx_{n-i}}x_{n-i}
 \label{eq:Mcgrath2006_generalized}
\end{equation}
since $z_n=\frac{x_n}{x_{n-i-1}}$ reduce (\ref{eq:Mcgrath2006_generalized}) to 
$$
z_{n+1}=\frac{a+bz_{n-i}}{c+dz_{n-i}}
$$
that can be linearized to an equation of order $2i+2$ as we do for (\ref{eq:BajoLiz2011_generalized2}) (open problem \ref{op:Mcgrath2006_generalized}).\\

The difference equation
\begin{equation}
x_{n+1}=\frac{x_{n-2}}{\pm 1+x_nx_{n-1}x_{n-2}}
 \label{eq:Khalaf_Allah}
\end{equation}
of reference \cite{2009_Khalaf_Allah} becomes $u_{n+1}=\frac{u_n}{\pm1+u_n}$ with the change $u_n=x_nx_{n-1}x_{n-2}$. Generalizing, in \cite{Shojaei2011} the RDE of order $k+1$
\begin{equation}
x_{n+1}=\frac{\alpha x_{n-k}}{\beta + \gamma x_n\ldots x_{n-k}},~~~\alpha,\beta,\gamma>0
 \label{eq:Shojaei2011}
\end{equation}
is reduced to the Riccati form
\begin{equation}
z_{n+1}=\frac{\alpha z_n}{\beta+\gamma z_n}
 \label{eq:Riccati_Shojaei}
\end{equation}
when we use the multiplicative change of variables $z_n=x_n\ldots x_{n-k}$.\\
Let's generalize the FS characterization of \cite{Shojaei2011} to the case where $\alpha,\beta,\gamma$ are arbitrary complex numbers. If one of them is zero, then equation (\ref{eq:Shojaei2011}) becomes trivial or globally periodic and $\mathcal F$ is empty or the set of points with at least one zero component.\\
If $\alpha\beta\gamma\neq0$, let $c=\frac{\alpha}{\beta}$. Therefore (\ref{eq:Riccati_Shojaei}) transforms into
$$z_{n+1}=\frac{cz_n}{1+z_n}$$
The former equation admits the following closed form solution
$$z_{n}=\frac{c^nz_0}{1+(1+c+\ldots+c^{n-1})z_0}$$
And from here it is easy to give the forbidden set expression. If $c=1$, then $z_n=\frac{z_0}{1+nz_0}$ and $\mathcal F=\{\frac{-1}n:n\geq1\}$. If $c$ is another root of the unity, then the equation is globally periodic and the forbidden set of the Riccati equation is finite:
$$\mathcal F=\left\{-1,\frac{1-c}{c^2-1},\ldots,\frac{1-c}{c^{m-1}-1}\right\}$$
supposing $m$ to be the smallest positive integer such that $c^m=1$.\\
When we apply these results to (\ref{eq:Shojaei2011}), we get that
%
%
%
\begin{equation}
\mathcal F=\bigcup_{n=0}^{M}\left\{(x_{-k},\ldots,x_0)\in\R^{k+1}:x_{-k}\ldots x_0=\frac{-1}{\sum_{i=0}^n c^i} \right\} 
\label{eq:FS_Shojaei2011}
\end{equation}
that is, a countable union of generalized hyperbolas when $c$ is not a root of the unity ($M=+\8$) or a finite union if $c^{M+1}=1$.\\

There are a number of RDEs for which the closed form solution and FS is given in \cite{rhouma2005closed}. Those equations can be grouped in three categories:
\begin{itemize}
 \item Those obtained from the multiplicative DE
  \begin{equation}
   x_{n+1}=ax_n^p x_{n-1}^q, ~~~ a\neq0~~~ p,q\in\Z
   \label{eq:Rhouma_multiplicative}
  \end{equation}
  when we do change of variables as 
  $$x_n=y_n-b$$
  $$x_n=\frac{y_n+b}{y_n+c}$$
  We deal then with RDE of order $2$ and degree $2$. It is not difficult to obtain the closed form of (\ref{eq:Rhouma_multiplicative}) from which the FS expression is constructed.
  
  \item RDE of order $2$ and degree $2$ resulting of the introduction of the variable $x_n=\frac1{y_n-k}$ in the general linear DE
  \begin{equation}
  x_{n+1}=ax_n+bx_{n-1}+c
  \label{eq:Rhouma_linear}
  \end{equation}
  \item And another family of RDE of order $2$ and degree $2$, given in this case by the change $x_n=\frac{y_{n+1}+\alpha y_n+\beta}{\gamma y_{n+1}+\lambda y_n+\mu}$ applied to the Riccati DE of order $1$
  \begin{equation}
  x_{n+1}=\frac{ax_n+b}{cx_{n}+d}
  \label{eq:Rhouma_Riccati}
  \end{equation}
\end{itemize}

In each of the former cases we propose some generalizations. Given a M\"obius transformation $T(x)=\frac{\alpha x +\beta}{\gamma x+\delta}$, it must be possible to explicitly determine the FS of the RDE constructed with the change $x_n=T(y_n)$ applied to equation (\ref{eq:Rhouma_multiplicative}) or applied to equation (\ref{eq:Rhouma_linear}). These are the claims of open problems \ref{op:Rhouma_multiplicative} and \ref{op:Rhouma_linear} respectively.

In the third case, we have that the closed form solution of (\ref{eq:Rhouma_Riccati}) and the change of variables able to write
$$x_n=g(a,b,c,d,y_0,y_{-1})=g_n$$

Moreover, the relationship $x_n=\frac{y_{n+1}+\alpha y_n+\beta}{\gamma y_{n+1}+\lambda y_n+\mu}$ can be interpreted as a nonautonomous linear DE when solved for $y_{n+1}$, that is
$$y_{n+1}=\frac{\lambda g_n-\alpha}{1-\gamma g_n}y_n+\frac{\mu g_n-\beta}{1-\gamma g_n}
$$

Therefore a general expression for $y_n$ and for $\mathcal F$ can be computed, and of course, the same idea must work for every change of variables $x_n=H(y_{n+1},\ldots,y_{n-k})$ where the explicit solution of the nonautonomous $y_n$ equation is known, and for every equation $x_{n+1}=G(x_n,\ldots,x_{n-k})$ whose closed form is also known (open problem \ref{op:Rhouma_generalized}).

\section{Use of invariants}
An interesting modification of the former ideas is the use of invariants to describe the FS of some RDEs. Consider the following example from \cite{Palladino2012}
\begin{equation}
x_{n+1}=\frac{x_n}{1+Bx_{n-1}-Bx_n}
 \label{eq:Palladino_1}
\end{equation}
where $B\in\C\setminus\{0\}$. This equation has the following invariant
\begin{equation}
\left(\frac 1 {x_n}+B\right)(1+Bx_{n-1})=C
\label{eq:Palladino_1_invariant}
\end{equation}
that is, for every solution $(x_n)_{n=-1}^{+\8}$ of (\ref{eq:Palladino_1}) there exist $C=C(x_{0},x_{-1})$ such that (\ref{eq:Palladino_1_invariant}) holds for every $n\geq0$.\\

The presence of an invariant ables to write an alternative form of the DE. Indeed, solving (\ref{eq:Palladino_1_invariant}) for $x_{n}$ and changing the indices, we get the following Riccati DE
$$x_{n+1}=\frac{1+Bx_n}{C-B-B^2x_n}$$
where recall that constant $C$ depend on the values of $x_{-1}$ and $x_0$. This is the basic idea to describe the FS of (\ref{eq:Palladino_1}).\\
Note that formula (\ref{eq:Palladino_1_invariant}) implies the identity
$$\left(\frac 1 {x_{n+1}}+B\right)(1+Bx_{n})=\left(\frac 1 {x_n}+B\right)(1+Bx_{n-1})$$
from where (\ref{eq:Palladino_1}) is deduced.\\
A generalization of the former remark is the following. Let $T_1$ and $T_2$ be two M\"obius transformations, $T_i(x)=\frac{\alpha_ix+\beta_i}{\gamma_ix+\delta_i}$, $i=1,2$, and consider the invariant
\begin{equation}
 T_1(x_n)T_2(x_{n-1})=C
 \label{eq:Palladino_generalized_invariant}
\end{equation}
Therefore 
$$T_1(x_{n+1})T_2(x_{n})=T_1(x_n)T_2(x_{n-1})$$
and from here we get
\begin{equation}
 x_{n+1}=T_1^{-1}\left(\frac{T_1(x_n)T_2(x_{n-1})}{T_2(x_n)}\right)
 \label{eq:Palladino_generalized_equation}
\end{equation}
a RDE of order $2$ whose FS could be described with the former methodology as the invariant (\ref{eq:Palladino_generalized_invariant}) implies also that
$$
x_{n+1}=T_1^{-1}\left(\frac{C}{T_2(x_{n})}\right)
$$
which is a Riccati DE of order $1$ (open problem \ref{op:Palladino_generalized}).

For example, given $B_1$ and $B_2$ non zero complex numbers, the invariant $ (x_{n}+B_1)(x_{n-1}+B_2)=C$
leads to DE
$$x_{n+1}=\frac{x_nx_{n-1}+(B_2-B_1)x_n+B_1x_{n-1}}{x_n+B_2}$$
and the invariant $ \frac{x_{n}+B_1}{x_{n-1}+1}=C$
to DE
$$x_{n+1}=\frac{x_n^2+(1+B_1)x_n-B_1x_{n-1}}{x_{n-1}+1}$$
both of them not included in reference \cite{Palladino2012}.\\

Also, given $k\geq1$, the invariant
\begin{equation}
 T_1(x_n)T_2(x_{n-k})=C
 \label{eq:Palladino_generalized_invariant_order_k}
\end{equation}
produces the DE
\begin{equation}
 x_{n+1}=T_1^{-1}\left(\frac{T_1(x_n)T_2(x_{n-k+1})}{T_2(x_n)}\right)
 \label{eq:Palladino_generalized_equation_order_k}
\end{equation}
that can be studied in the same way because the $C$ dependent DE
$$
x_{n+1}=T_1^{-1}\left(\frac{C}{T_2(x_{n-k+1})}\right)
$$
is reduced to a linear equation of order $k+1$ using the same changes as in equation (\ref{eq:BajoLiz2011_generalized2}).\\

Another example of use of invariants is in \cite{2010_AS_OnTheAsymptotic}. DE
\begin{equation}
x_{n+1}=\frac{x_n^2+x_{n-1}^2-x_n(x_{n-1}+x_{n-2})}{x_{n-1}-x_{n-2}} 
\label{eq:AghajaniShouli}
\end{equation}
verifyes that $I(x_n,x_{n-1},x_{n-1})=\frac{x_n-x_{n-2}}{x_n-x_{n-1}}$ is constant over each solution of (\ref{eq:AghajaniShouli}). Therefore, from the equality
\begin{equation}
 \frac{x_{n+1}-x_{n-1}}{x_{n+1}-x_{n}}=C
\end{equation}
the following linear relation is deduced
\begin{equation}
 x_{n+1}=\frac{C}{C-1}x_n- \frac 1{C-1}x_{n-1}
\end{equation}
and from here it is easy to give the closed and forbidden set of (\ref{eq:AghajaniShouli}). In particular, the FS is
$$\mathcal F=\{(x,y,z)\in\R^3:x=y\text{ or }y=z\}$$
A possible generalization is to consider invariants of the form
\begin{equation}
  \frac{x_{n+1}-x_{n-k}}{x_{n+1}-x_{n-l}}=C~~~~k,l\in\N,~k\neq l
\label{eq:AghajaniShouli_invariant}
\end{equation}
as every solution will be associated to a linear recurrence (open problem \ref{op:AghajaniShouli}).\\

Finally, a general question is to determine which relationship exists between the poles and zeros of an algebraic invariant, from one side, and the elements of its associated DE, from other.
\section{Forbidden set curves given in an explicit form}
\label{sec:FS_ExplicitForm}
One of the oldest examples in the literature concerning the forbidden set problem is in \cite{TheForbiddenSet_CDV}. Let $p\leq-1$ be a real number, and consider the following RDE of order $2$
\begin{equation}
x_{n+1}=p + \frac{x_{n-1}}{x_n}
 \label{eq:CamouzisDeVault}
\end{equation}
Map $F(x,y)=\left(y,p+\frac x y\right)$ is the unfolding of (\ref{eq:CamouzisDeVault}) as we have that $F(x_{n-1},x_n)=(x_n,x_{n+1})$, $\forall n\geq0$. Let $G(x,y)=(x(y-p),x)$ be its inverse map, and $A=\{(x,0):x\in\R\}$. Therefore
$$\mathcal F=\bigcup_{n=0}^{+\8} G^n(A)$$
This is an obvious characterization of the FS as the set of inverse orbits of poles of the iteration map that can be improved as follows. Let's consider the subset $A^+=\{(x,0):x\geq0\}$ and let $(g_n)_{n=1}^{+\8}$ be the sequence of functions $g_n:[0,+\8)\rightarrow \R$ defined inductively using the operator $h_g(x)=x\cdot(-p+g^{-1}(x))$ in the following way
$$g_1(x)=-px, ~~~~ g_{n+1}(x)=h_{g_n}(x)$$
Therefore we have
\begin{teo}[\cite{TheForbiddenSet_CDV} theorem 2]
\label{th:TheForbiddenSet_CDV}
Let $\mathcal F^+$ be the subset of the forbidden set of (\ref{eq:CamouzisDeVault}) defined by $\mathcal F^+=\bigcup_{n=0}^{+\8} G^n(A^+)$. Then
\begin{equation}
 \mathcal F^+=A^+\cup\{(x,y)\in[0,+\8)\times\R: \exists n\geq1, y=g_n^{-1}(x)\}
 \end{equation}
 \end{teo}

\begin{figure}[h]
\begin{center}$
\begin{array}{cc}
\includegraphics[width=0.45\textwidth]{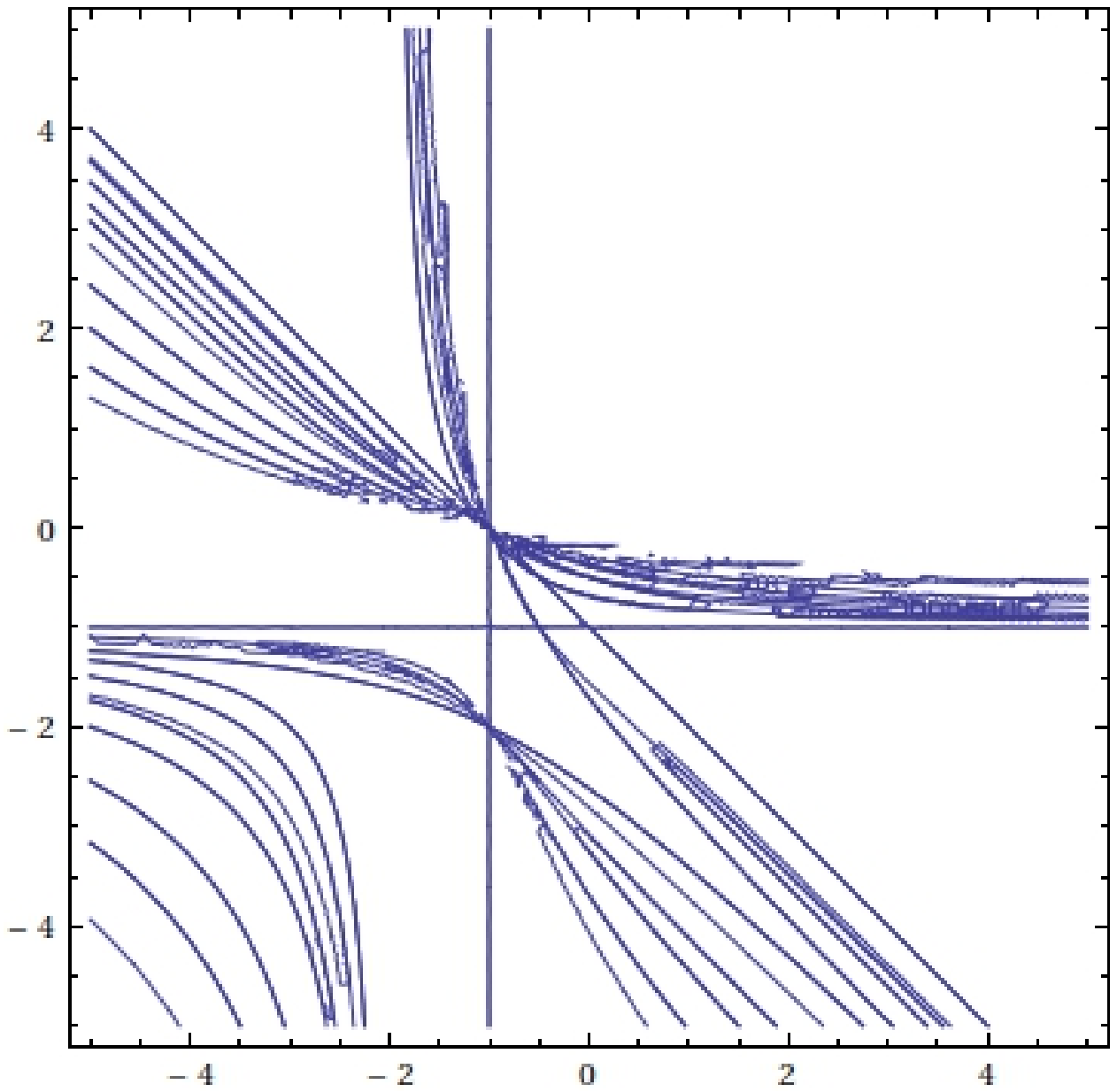} &
\includegraphics[width=0.45\textwidth]{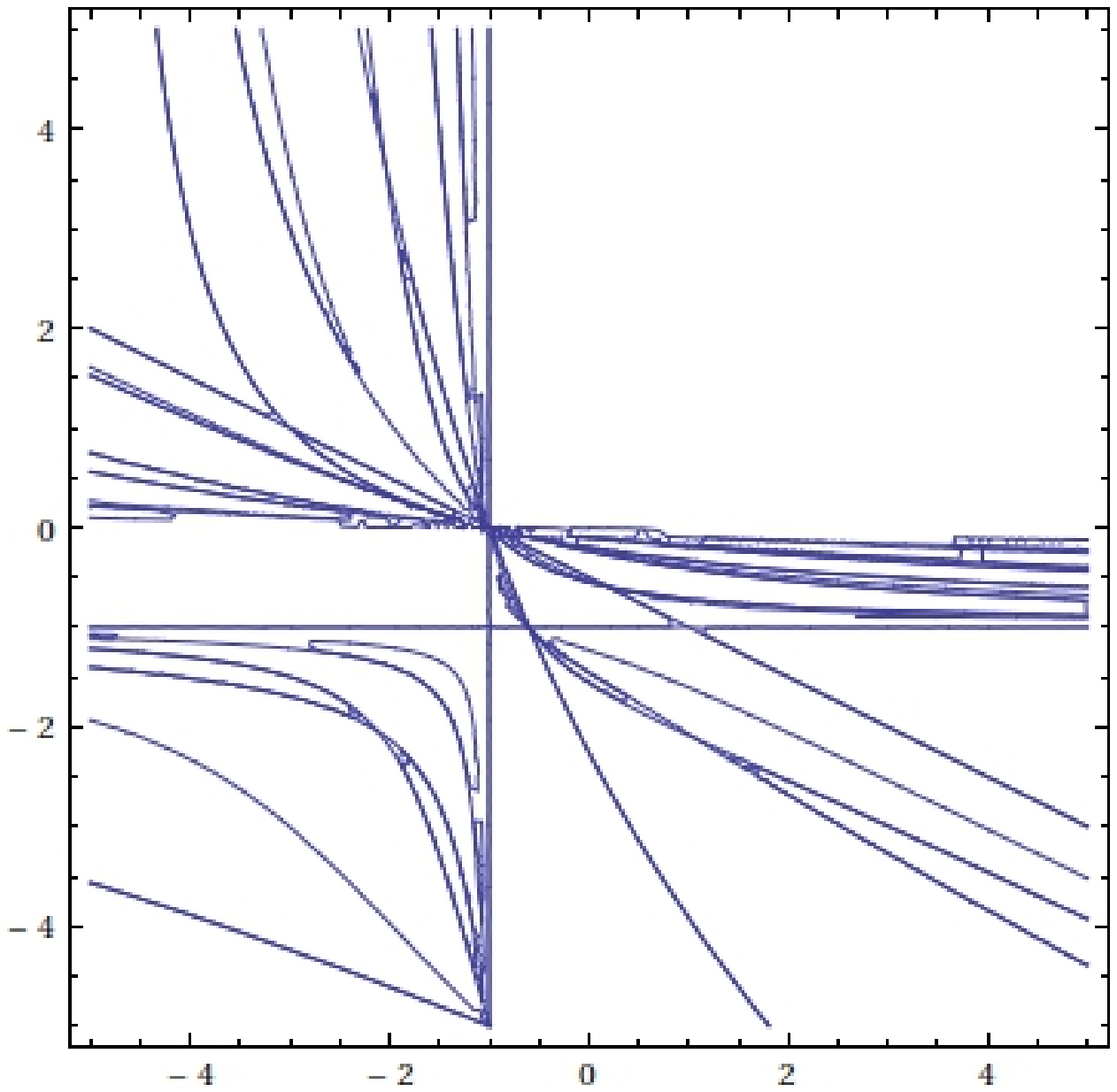}\\
\includegraphics[width=0.45\textwidth]{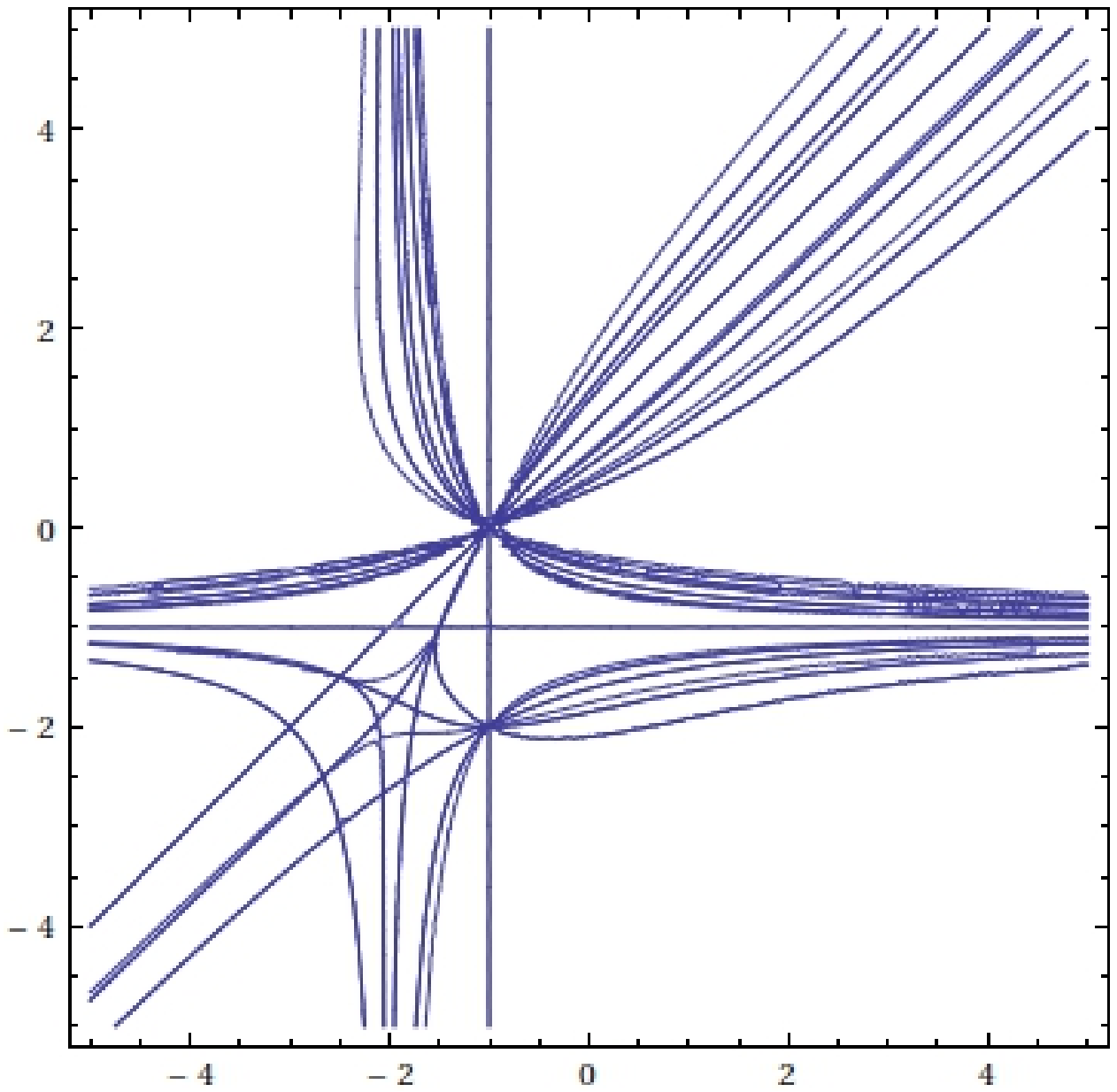} &
\includegraphics[width=0.45\textwidth]{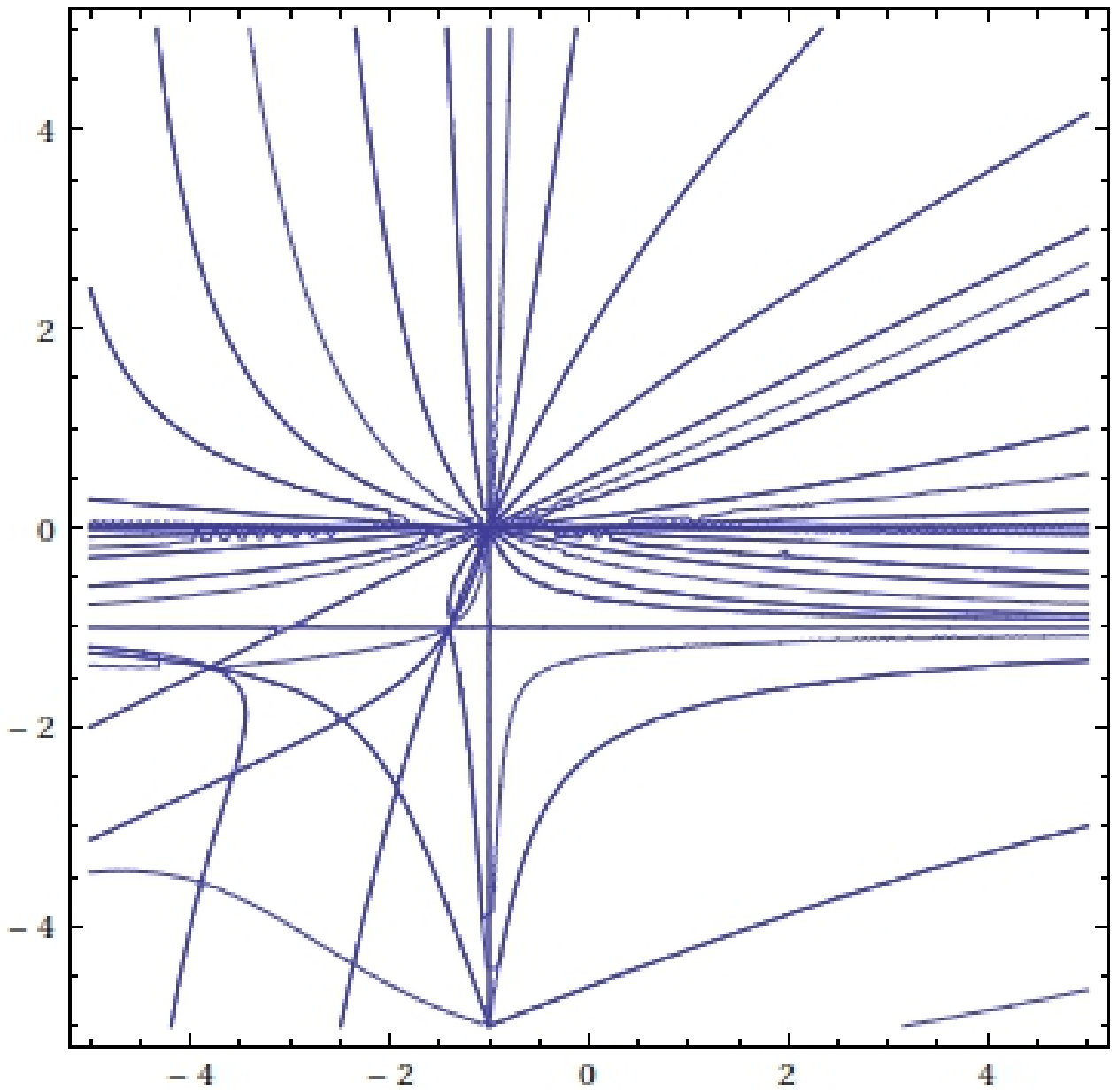}
\end{array}$
\end{center}
\caption{The first $10$ forbidden curves in Pielou's RDE. From up to bottom and left to right, $a=1$, $a=2$, $a=-1$ and $a=-2$}
\label{fg:Pielou}
\end{figure}

In general, given a RDE we can describe its FS as a set of implicitly defined hypersurfaces. In the second order case we get a set of implicitly defined curves. For example, in the case of Pielou's equation
\begin{equation}
x_{n+1}=\frac{a x_n}{1+x_{n-1}}
 \label{eq:Pielou}
\end{equation}
we can consider the unfolding map $F(x,y)=\left(y,\frac{ay}{1+x}\right)$ such that $F(x_{n-1},x_n)=(x_n,x_{n+1})$. By iterating map $F$ we get
$$F^2(x,y)=\left(\frac{a y}{1+x},\frac{a^2 y}{(1+x) (1+y)}\right)$$
$$F^3(x,y)=\left(\frac{a^2 y}{(1+x) (1+y)},\frac{a^3 y}{(1+y) (1+x+a y)}\right)$$
$$F^4(x,y)=\left(\frac{a^3 y}{(1+y) (1+x+a y)},\frac{a^4 (1+x) y}{(1+x+a y) \left(1+x+y+a^2 y+x y\right)}\right)$$
$$\vdots$$
From where the following first forbidden curves are deduced (see figure \ref{fg:Pielou}):
$$\begin{array}{lll}
\mathcal F &=&	\{(x,y)\in\R^2:1+x=0 \}~~\cup\\
	  & &	\{(x,y)\in\R^2:1+y=0 \}~~\cup\\
	  & &	\{(x,y)\in\R^2:1+x+a y=0 \}~~\cup\\
	  & &	\{(x,y)\in\R^2:1+x+y+a^2 y+x y=0 \}~~\cup\\
	  & &	\ldots
\end{array}$$
Remark that those curves are given in an implicit form. We propose to find a way to construct the explicit form as in theorem \ref{th:TheForbiddenSet_CDV} (open problem \ref{op:ExplicitCurves}).
\section{Symbolic description of the FS}
\label{sec:Symbolic}
\begin{figure}[h]
\begin{center}
\includegraphics[width=0.6\textwidth]{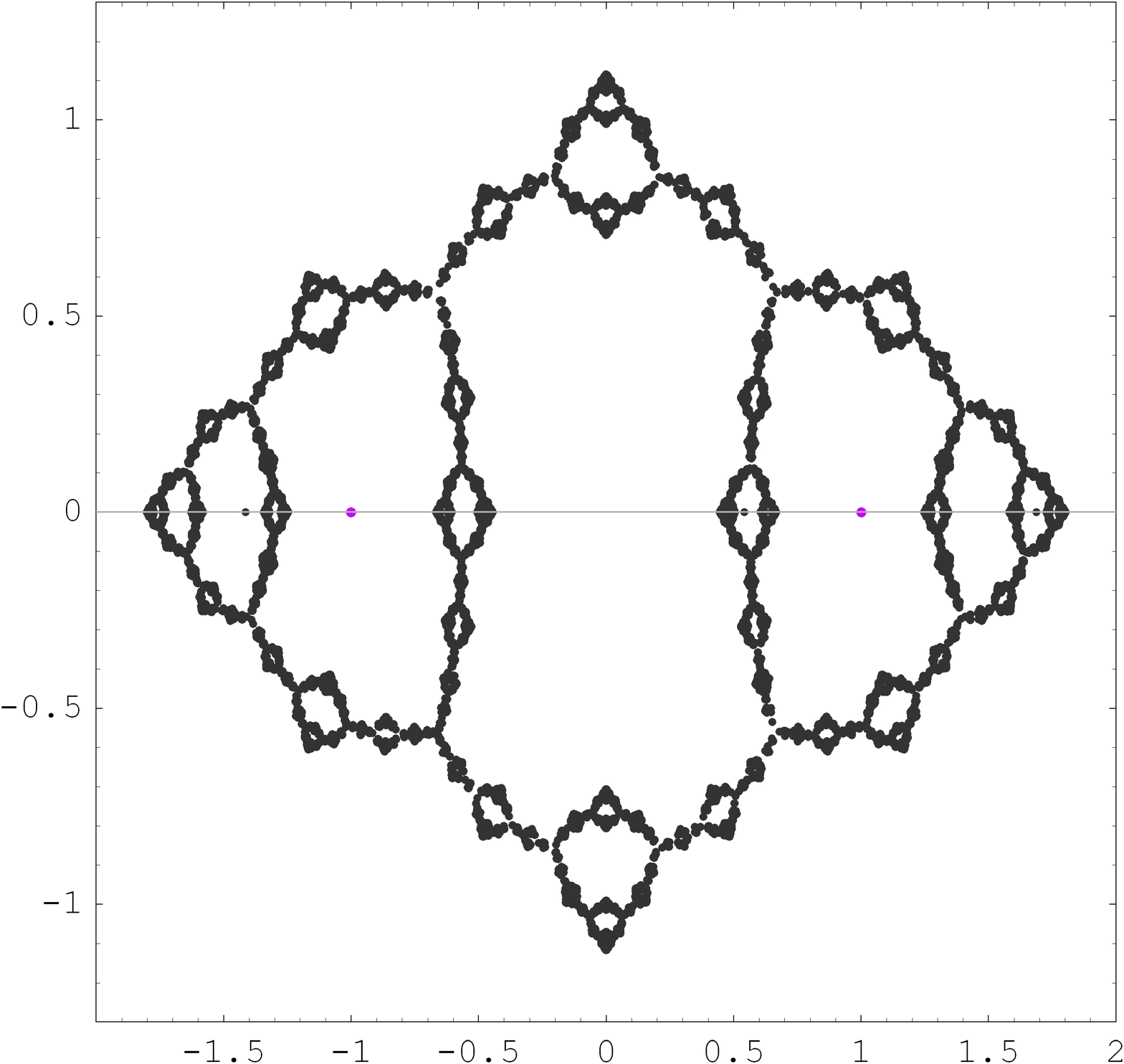}\\
\end{center}
\caption{Forbidden set of $x_{k+1}=\frac1{x_k^2-1}$.}
\label{fg:ForbiddenSetInverseParabola_C}
\end{figure}

Let's consider the following RDEs of order $1$ and degree $2$
\begin{equation}
x_{n+1}=\frac1{x_n^2-a}
 \label{eq:InverseParabola}
\end{equation}
\begin{equation}
x_{n+1}=\frac1{ax_n(1-x_n)}
 \label{eq:InverseLogistics}
\end{equation}
where $a$ is a real or complex parameter. These equations are studied in \cite{tesisACV} in connection with their relationship with Li-Yorke chaos in RDEs.\\
It is possible to describe the forbidden set of the particular case
\begin{equation}
x_{n+1}=\frac1{x_n^2-1}
 \label{eq:InverseParabola_a1}
\end{equation}
using a symbolic notation. Let 
\begin{equation}
\left\{
 \begin{array}{lll}
 h_+(x)&=&\sqrt{\frac1x+1}\\
h_-(x)&=&-\sqrt{\frac1x+1} 
 \end{array}
\right.
\label{eq:InverseBranches}
\end{equation}
be the inverse branches of the iteration function in (\ref{eq:InverseParabola_a1}). Clearly, $\mathcal F$ is the forward multiorbit of the poles of $f(x)=\frac1{x^2-1}$. Moreover
\begin{pro}
The unique preimage of pole $-1$ in (\ref{eq:InverseParabola_a1}) is $0$. The preimages of pole $1$ are $a_1\ldots a_n$,
where $a_i=h_+$ or $a_i=h_-$ for $i=1,\ldots,n$, and $a_1\ldots a_n$ is the abbreviation of $a_1\circ \ldots \circ a_n(1)$. Those preimages are always well defined in $\C$, that is, $a_1\ldots a_n\neq 0$ for every word, and the former representation is unique, that is, given two words such that $a_1\ldots a_n=b_1\ldots b_m$ then $n=m$ and $a_i=b_i$ for $i=1,\ldots,n$.\\
Therefore
\begin{equation} 
\mathcal F=\{-1,0,1\}\cup\{a_1\ldots a_n: n\in\N, a_i=h_{\pm}\}\subset\C
 \label{eq:ForbiddenSetInverseParabola_C}
\end{equation}
\end{pro}
\begin{dem}
Let's see that words $a_1\ldots a_n$ are distint from zero and are unique, as the remaining assertions are straightforward to demonstrate.\\
When $n=1$, $a_1(1)=\pm\sqrt2\neq0$. For $n>1$, if $a_1\ldots a_n=0$, then $a_1(a_2\ldots a_n)=0$ and therefore $a_2\ldots a_n=-1$, what is impossible as $-1$ does not belong to the image of $h_+$ or $h_-$.\\
Now, suppose that $a_1\ldots a_n=b_1\ldots b_m$. There is no lack of generality if we assume that $a_1\neq b_1$, because $h_+$ and $h_-$ are injective maps. Let $x=a_2\ldots a_n$ and $y=b_2\ldots b_m$. The former equality is rewritten as $a_1(x)=b_1(y)$. But note that $\Im h_+\cap \Im h_-=\{0\}$ from where we deduce that $x=y=-1$, what is again impossible. To avoid the contradictions both words $a_1\ldots a_n$ and $b_1\ldots b_m$ must be equal.
\end{dem}
Note that this is a set in the complex field (see figure \ref{fg:ForbiddenSetInverseParabola_C}). If we center in the real field, note that applying $h_-$ one time we stay over the reals, and applying it two consecutive times we get $\C\setminus\R$. Therefore
\begin{cor}
 The forbidden set of (\ref{eq:InverseParabola_a1}) over $\R$ is
$$\mathcal F_{\R}=\{-1,0,1\}\cup\{a_1\ldots a_n: n\in\N, a_i=h_{\pm}, \text{ there is no two consecutive } h_-\}$$
\end{cor}
Numerically, $\mathcal F_{\R}$ is the projection of figure \ref{fg:ForbiddenSetInverseParabola_C} over the real line.
In open problem \ref{op:SymbolicDescription} we propose to find the symbolic description of $\mathcal F_{\C}$ and $\mathcal F_{\R}$ for the families of RDEs (\ref{eq:InverseParabola}) and (\ref{eq:InverseLogistics}).
 
\section{Qualitative approaches}
\label{sec:Qualitative}
As we have seen in the former sections, there are some works describing the FSs of DEs. But it should be remarked that the most part of DEs don't have a closed form and therefore there is probably no explicit expression for their FSs.\\
We claim that new perspectives should be taken in order to solve the problem of the FS for a generic DE. Specifically, given a DE we propose three general questions about its FS
\begin{enumerate}
 \item Which are its topological and metrical properties? (Openness, Boundedness, Fractal dimension,...)
 \item What can be said about perturbations of its parameters?\\
 For example, given a family of DEs $x_{n+1}=f_a(x_n)$ depending on parameter $a$, and supposing that some feature of the FS is known when $a=a_0$, what can be said about the FS of the DE when $a\in(a_0-\varepsilon,a_0+\varepsilon)$?
 \item Are there some relationships between the properties of the FS and the dynamical behavior of the DE?
\end{enumerate}
There are, as far as we know, very few papers with a qualitative approach. Let's briefly review them.\\

\begin{figure}[h]
\begin{center}
\includegraphics[width=0.4\textwidth]{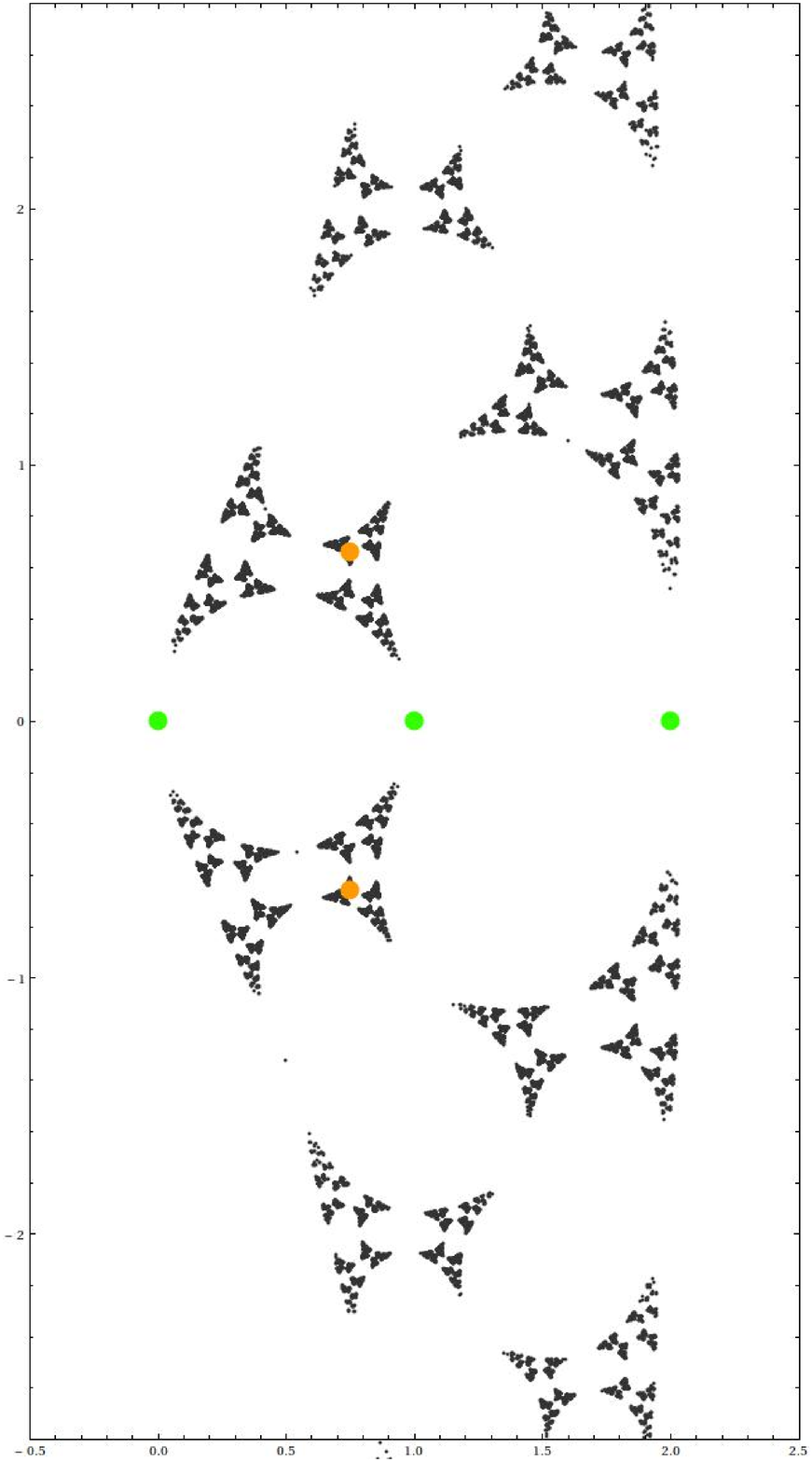}\\
\end{center}
\caption{$50000$ elements of the FS of (\ref{eq:RubioMassegu_CE}). Green points are the unique elements of $\mathcal F$ in $\R$. Fixed points are marked in yellow, and they seem to be in the clousure of $\mathcal F$.}
\label{fg:FS_RubioMassegu}
\end{figure}

In \cite{RubioMassegu2009} the problem of closedness of the FS is discussed in connection with the global periodicity question. Let $\mathbb K$ be the real or complex field. Let $F:\mathbb K^k\rightarrow\mathbb K^k$ and consider the $k$-dimensional DE or dynamical system
\begin{equation}
x_{n+1}=F(x_n)
 \label{eq:RubioMassegu_Map}
\end{equation}
Let $D\subset \mathbb K^k$ be the domain of $F$. A set $A\subset D$ has full Lebesgue measure in $D$ if the measure of $D\setminus A$ is zero. The map $F$ is almost a local diffeomorphism if there exists an open set $V\subset \mathbb K^k$ of full Lebesgue measure in $D$ such that $F$ is continuously differentiable on $V$ (or analytic on $V$ when $\mathbb K = \C$) and $DF(x)\neq0$ a.e. $x\in V$.
\begin{teo}[\cite{RubioMassegu2009}, theorem 3.]
Let $F:D\subseteq \mathbb K^k\rightarrow \mathbb K^k$ be an almost local diffeomorphism such that the set $F^{-1}(\mathbb K^k\setminus D)$ has zero Lebesgue measure. Suppose also that $D$ is open and that $F$ is a continuous globally periodic map. Then the good set of (\ref{eq:RubioMassegu_Map}) is closed and has full measure in $D$.
 \label{teo:RubioMassegu_MainTheorem}
\end{teo}
There are some consequences of this result. The natural domain of a rational map $F=(F_1,\ldots,F_k)$ is the set of points in $\mathbb K^k$ such that any denominator in $F$ vanish. The first result is for the rational case.
\begin{cor}[\cite{RubioMassegu2009}, corollary 4]
Let $D$ be the natural domain of a rational map $F:D\subseteq \mathbb K^k\rightarrow \mathbb K^k$. If $F$ is globally periodic and $det DF(x)$ is not identically zero, then the good set of (\ref{eq:RubioMassegu_Map}) is open and has full Lebesgue measure in $\mathbb K^k$.
\end{cor}
In the case of RDE, we say that a function $f(x_1,\ldots,x_{k})$ depends effectively on its last variable if $\frac{\partial f}{\partial x_k}(x)\neq0$ for some $x=(x_1,\ldots,x_k)\in D$. Therefore for RDE
\begin{equation}
x_{n+1}=f(x_n,\ldots,x_{n-k+1})
 \label{eq:RubioMassegu_DE}
\end{equation}
we get
\begin{cor}[\cite{RubioMassegu2009}, corollary 5]
 Let $f:D\subseteq\mathbb K^k\rightarrow\mathbb K$ be a rational function, and $D$ its natural domain. If $f$ depends effectively on its last variable and (\ref{eq:RubioMassegu_DE}) is globally periodic, then the good set is open and has full Lebesgue measure in $\mathbb K^k$.
\end{cor}
\noindent A more concrete application is given in the case of the Riccati DE of order $1$
\begin{pro}[\cite{RubioMassegu2009}, proposition 7.]
DE (\ref{eq:RiccatiOrder1}) such that none of the quantities $d$, $ad-cb$ and $b+c$ is zero, has an open good set if, and only if, the equation is globally periodic.
\end{pro}
It is a natural question to wonder if this result applies also to any RDE. Again in \cite{RubioMassegu2009} is given the following not globally periodic example
\begin{equation}
 x_{n+1}=-\frac{(x_n-1)(x_n-2)}{x_n}
 \label{eq:RubioMassegu_CE}
\end{equation}
whose FS is $\mathcal F=\{0,1,2\}$ and therefore closed. Numerical experiments show also that $\mathcal F$ is not closed when the equation is regarded in $\C$ (see figure \ref{fg:FS_RubioMassegu}, open problem \ref{op:RubioMassegu} and conjecture \ref{conj:RubioMassegu}).\\

A second example of qualitative study of DE is in \cite{RubioMassegu2007OnTheExistence}. There the author stands sufficient conditions for a DE to have a good set of full Lebesgue measure, generalizing a previous remark of \cite[section 2]{2001_CsornyeiLaczkovich_SomePeriodic}. In the nonautonomous case, the main result is
\begin{teo}[\cite{RubioMassegu2007OnTheExistence}, theorem 5]
Let $f:D\subseteq\mathbb K^k\rightarrow\mathbb K$ be a function locally in fact of order $k$ and such that the set
$M=\{x\in D:(f(x),x_2,\ldots,x_{k-1})\}$ has Lebesgue measure zero. Then the good set of equation (\ref{eq:RubioMassegu_DE}) has full measure in $D$.
\end{teo}

Here a function $f$ is locally in fact of order $k$ if there exists an open set $V\subseteq D$ of full Lebesgue measure in $D$ and such that $\frac{\partial f}{\partial x_k}(x)\neq0$ a. e. $x\in V$.\\
In \cite{RubioMassegu2007OnTheExistence} it is remarked the importance of determining the topological and geometrical properties of the good set in terms of the iteration function and its domain, and two related open problems are proposed (see open problems \ref{op:RubioMassegu_EquilibriumInTheInterior} and \ref{op:RubioMassegu_OpenGoodSet}).\\

Also in \cite{Szalkai} there is a generalization of these considerations, studying the problem of avoiding forbidden sequences.\\

In \cite{Sedaghat2000} the following DE are considered:
\begin{equation}
x_{n+1}=\frac a{x_n^p}+1
\label{eq:Sedaghat_Qualitative_1}
\end{equation}
\begin{equation}
 x_{n+1}=\frac{a}{x_{n-1}^q}-\frac{1}{x_n^p}
\label{eq:Sedaghat_Qualitative_2}
\end{equation}
where $a\neq0$ and $p$ and $q$ are odd rationals with at least one of them positive. An odd rational is a fraction of the form $\frac{2i-1}{2j-1}$, $i,j\in\N$.\\

Again the objective is to determine somehow the forbidden set. The author uses the term \emph{crash set} referring to $\mathcal F$, emphasizing that the discontinuity is a pole.\\

\begin{figure}[h]
\begin{center}$
\begin{array}{cc}
\includegraphics[width=0.45\textwidth]{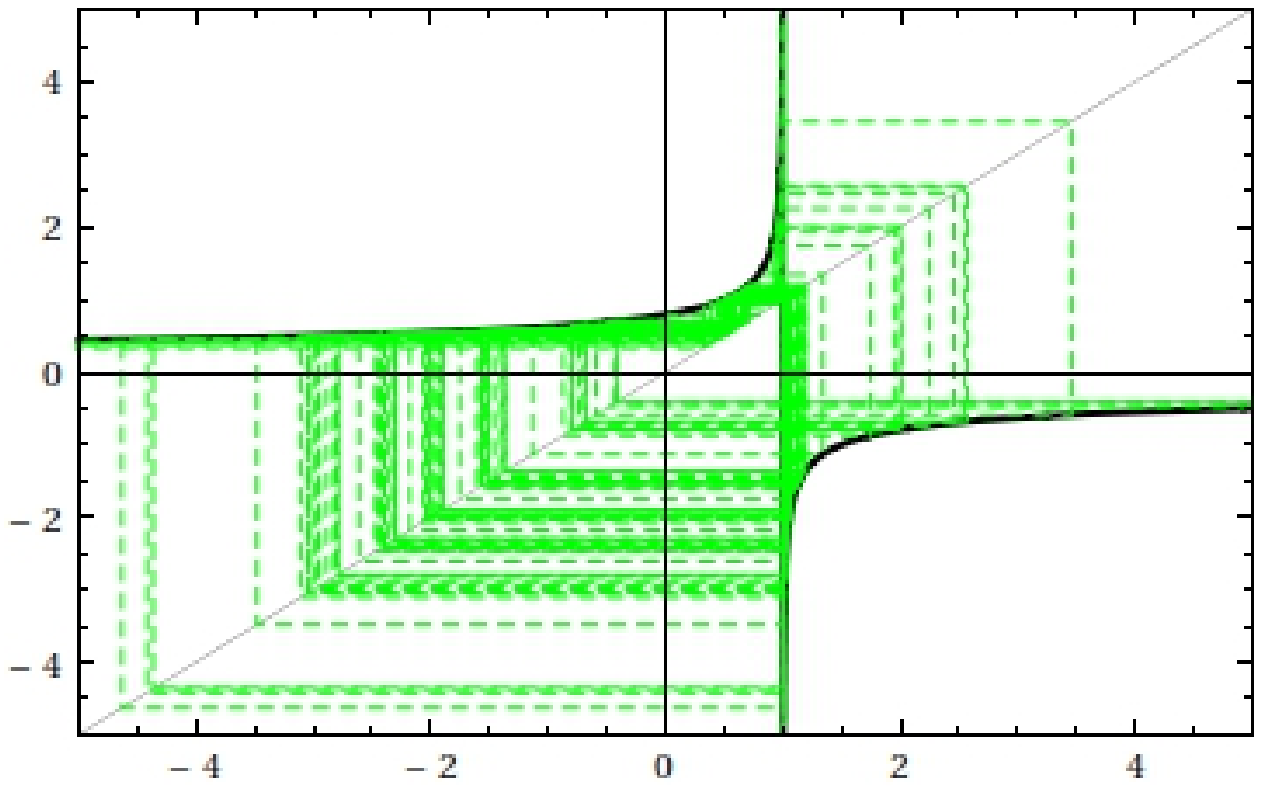} &
\includegraphics[width=0.45\textwidth]{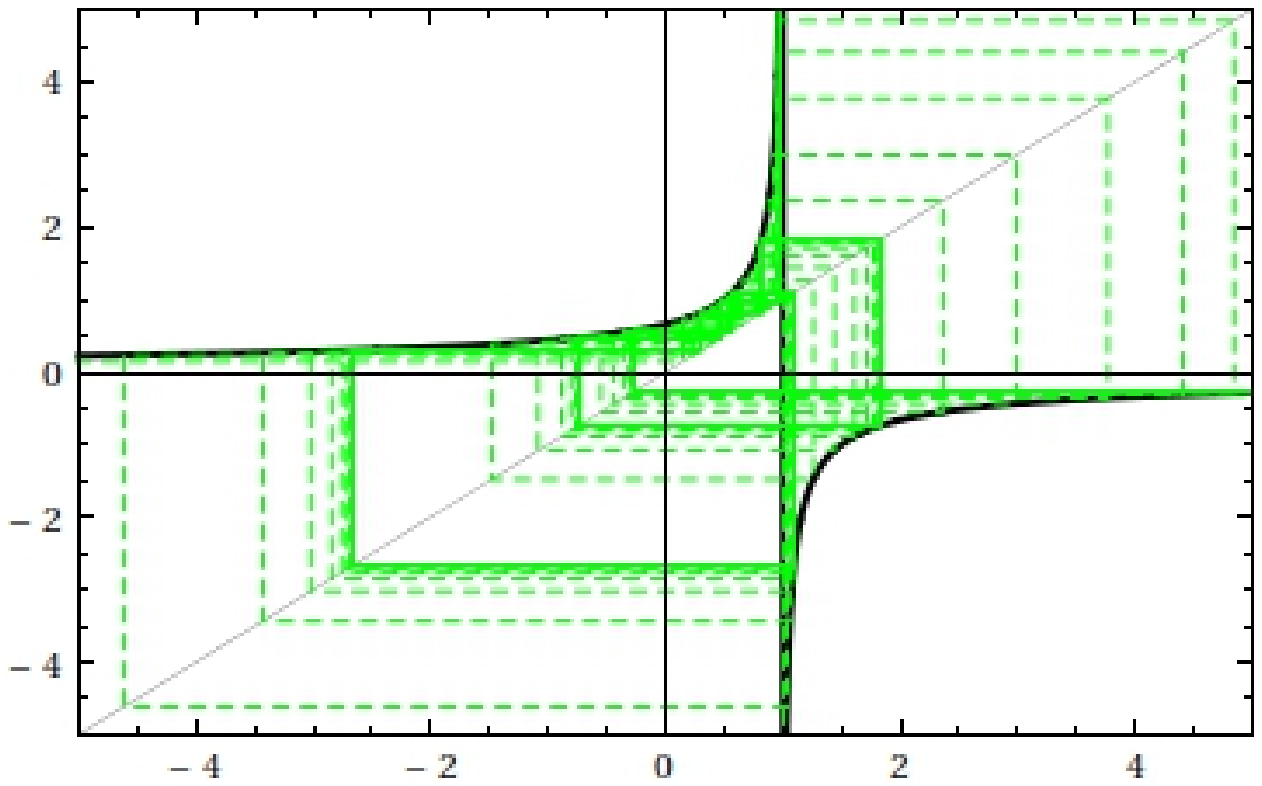}\\
\includegraphics[width=0.45\textwidth]{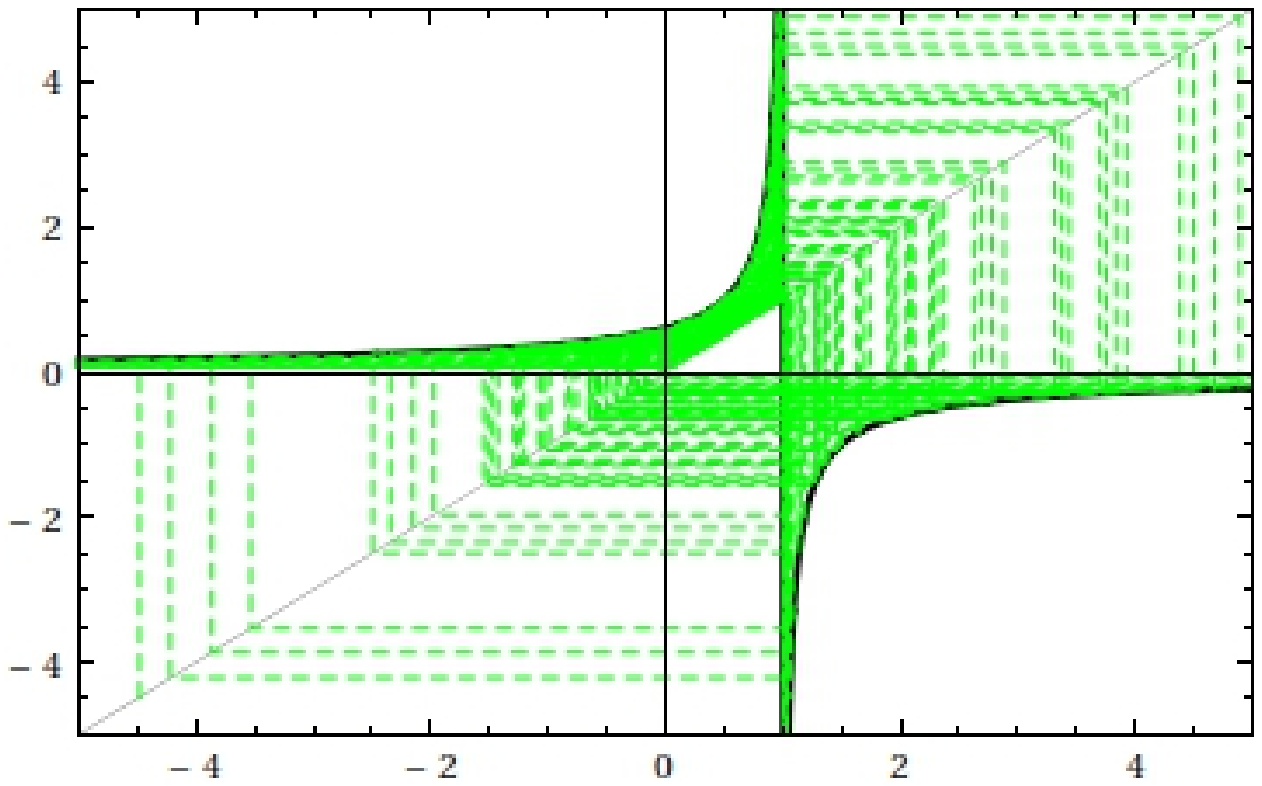} &
\includegraphics[width=0.45\textwidth]{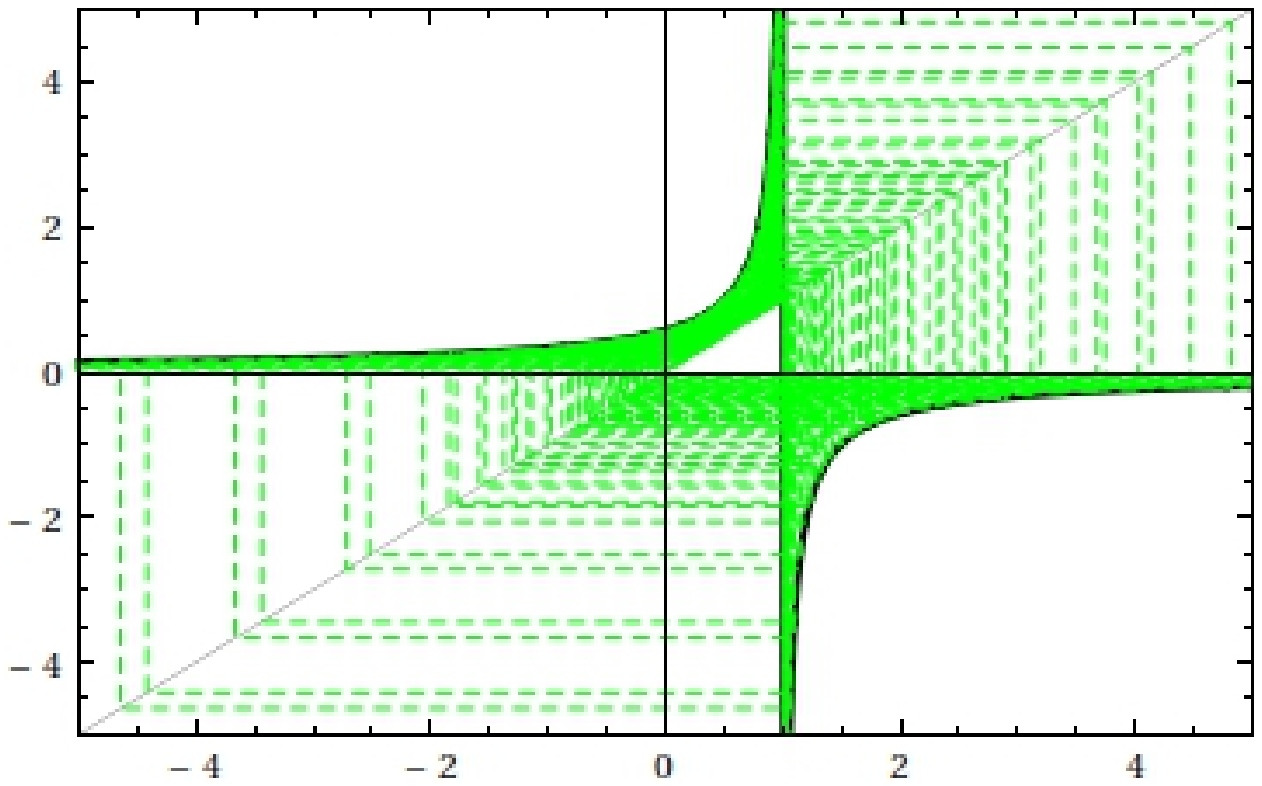}
\end{array}$
\end{center}
\caption{Forbidden sets of function $f(x)=1-\frac{1}{2x^p}$ represented with a cobweb diagram of the inverse iteration starting at the pole $0$. From left to right and up to bottom, $p$ equals to $3$, $\frac 5 3$, $\frac 7 5$ and $\frac 9 7$ respectively. $250$ iterations were taken in each case.}
\label{fg:cobweb}
\end{figure}

None of the forbidden or crash sets of the former equations is explictly computed, but it is described qualitatively. In the case of (\ref{eq:Sedaghat_Qualitative_1}), the iteration function belongs to the family of monotonic maps with pole, defined as the class of functions $f:\R\rightarrow\R$ verifying
\begin{enumerate}
 \item[(A1)] $f$ is continuous on its domain $\R\setminus\{0\}$
 \item[(A2)] $f$ is injective on $\R\setminus\{0\}$
 \item[(A3)] $0\in f(\R\setminus\{0\})$
 \item[(A4)] $\lim\limits_{x\rightarrow 0}|f(x)|=\8$
\end{enumerate}
Note that $f$ must be an increasing or a decreasing map, and that $\mathcal F$ consists in the collection of preimages of the pole $0$, qualitatively described in the next theorem. $P$ stands for the set of fixed points of $f$.
\begin{teo}[\cite{Sedaghat2000}, theorem 1]
Let $f$ be a monotonic map with pole.
 \begin{enumerate}[(a)]
  \item If $f$ is a decreasing map, let $\overline x$ be the fixed point in $(f^{-1}(0),0)$. Then the sequence $\mathcal F=(f^{-n}(0))_n$ is contained in the interval $[f^{-1}(0),f^{-2}(0)]$, and converges either to $\overline x$ or to a $2$-cycle.
  \item Suppose that $f$ is an increasing map, $P\neq\emptyset$, and $\overline x=\inf P>f^{-1}(0)$. Then $\mathcal F$ is contained in $[f^{-1}(0),\overline x)$ and converges monotonically to $\overline x$.
 \end{enumerate}
 \label{th:Sedaghat_MonotonicMapWithPole}
\end{teo}
Remark that there are some cases not included in this result. When $f$ is increasing and $\overline x=\sup P<f^{-1}(0)$, cobweb diagrams show that $\mathcal F$ is a subset of $(\overline x,f^{-1}(0)]$ and converges monotonically to $\overline x$. Whereas for $P=\emptyset$, $\mathcal F$ could be finite or could have a more complicated structure. For example in \cite{Sedaghat2000} is wondered if in the last case $\mathcal F$ can be dense outside a compact set. (See open problem \ref{op:Sedaghat_DenseFS} and figure \ref{fg:cobweb}).\\

Theorem \ref{th:Sedaghat_MonotonicMapWithPole} applies not only to equation (\ref{eq:Sedaghat_Qualitative_1}), but also to equations as:\\
\begin{equation}
x_{n+1}=\frac{a}{\sinh x_n}+1
 \label{eq:sinh}
\end{equation}
and, in general to DEs of the form
\begin{equation}
x_{n+1}=\frac{a}{\phi( x_n)}+1
 \label{eq:bijection}
\end{equation}
where $a\neq0$ and $\phi:\R\rightarrow\R$ is a bijection such that $\phi(0)=0$.\\

The FS of (\ref{eq:Sedaghat_Qualitative_2}) is given in \cite{Sedaghat2000} in terms of the values of $a$ and $p$. We propose to make a similar study for the family of DE
\begin{equation}
x_{n+1}=T(x_n^p)
 \label{eq:Sedaghat_Qualitative_1_generalized}
\end{equation}
where $p$ is an odd rational and $T$ is the M\"obius transformation $T(x)=\frac{ax+b}{cx}$ with $a,b,c\neq0$ (open problem \ref{op:Sedaghat_Qualitative_1_generalized}).\\

In the case of equation (\ref{eq:Sedaghat_Qualitative_2}), in \cite{Sedaghat2000} there is a description of the FS curves in an explicit form using a recurrent algorithm, in a similar way as it was made in section \ref{sec:FS_ExplicitForm}.\\

Additional examples of estimative approaches where the FS is graphically represented can be found in \cite{CMAarticle,JDEAarticle,tesisACV} for the following DEs
\begin{equation}
x_{n+1}=\frac1{x_n+x_{n-2}}
 \label{eq:JDEAarticle}
\end{equation}
\begin{equation}
x_{n+1}=\frac{x_{n-1}}{1+x_n}
 \label{eq:tesis_orden2}
\end{equation}

\begin{figure}[h]
\begin{center}$
\begin{array}{cc}
\includegraphics[width=0.45\textwidth]{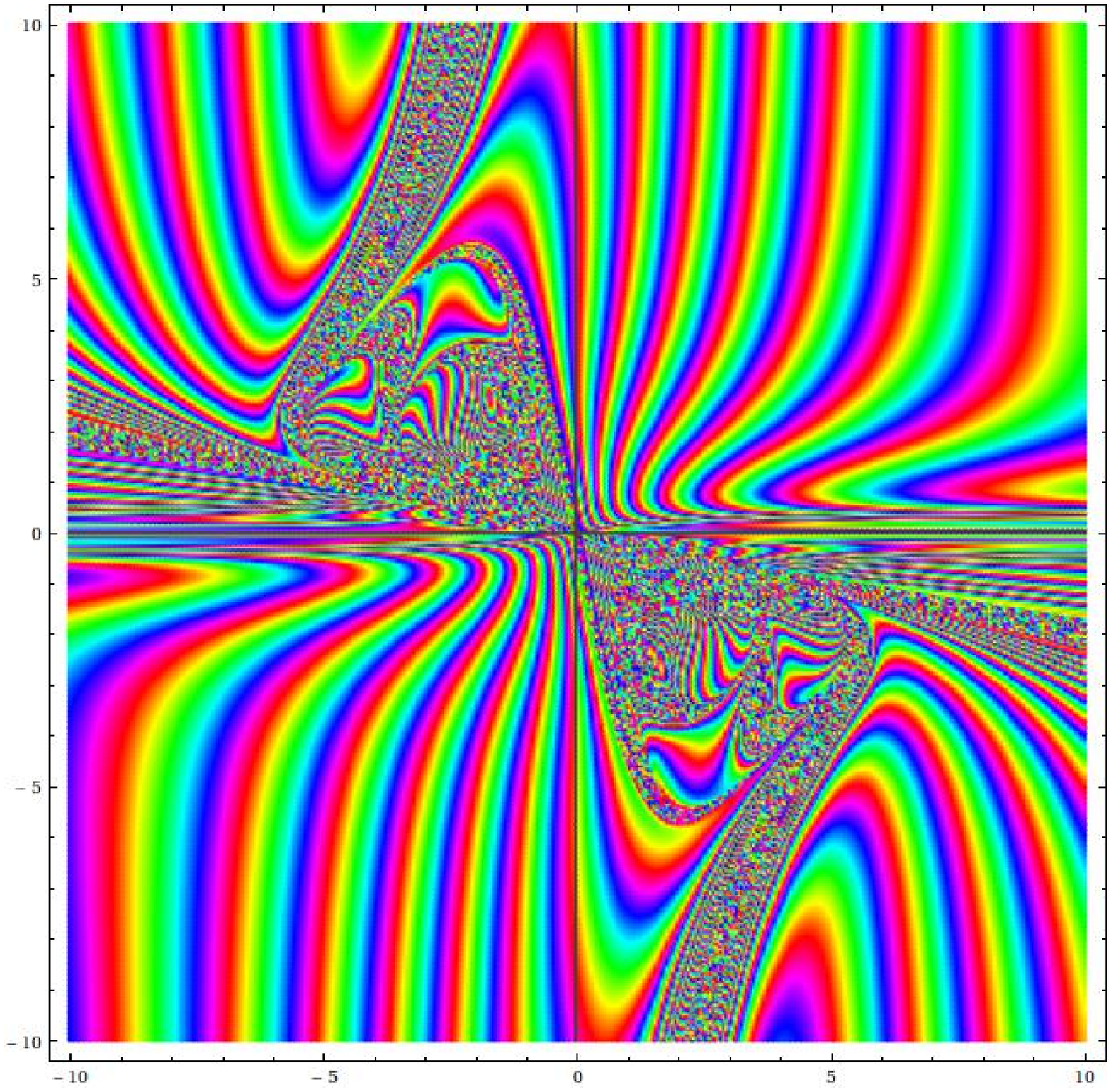} &
\includegraphics[width=0.45\textwidth]{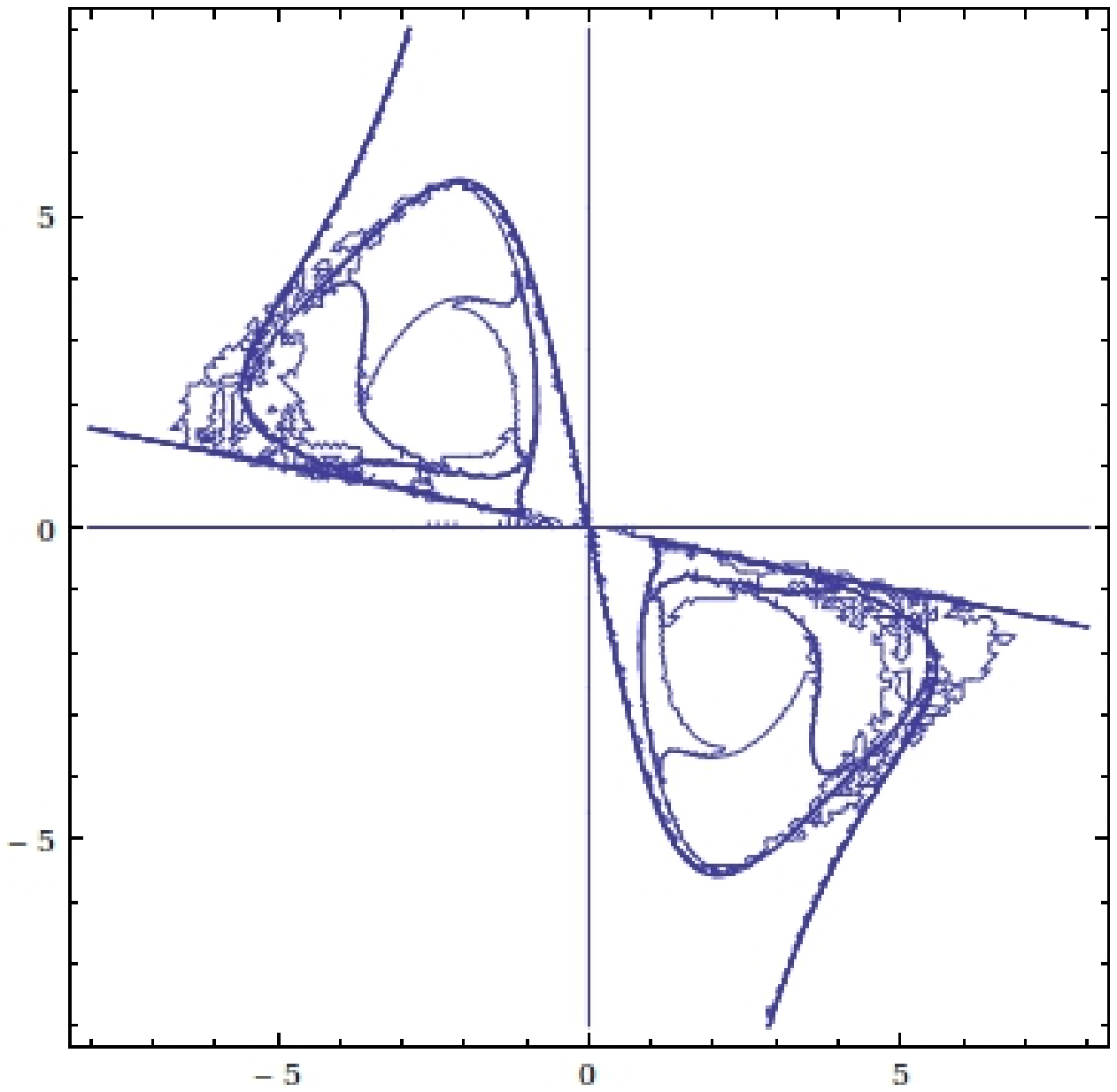}\\
\includegraphics[width=0.45\textwidth]{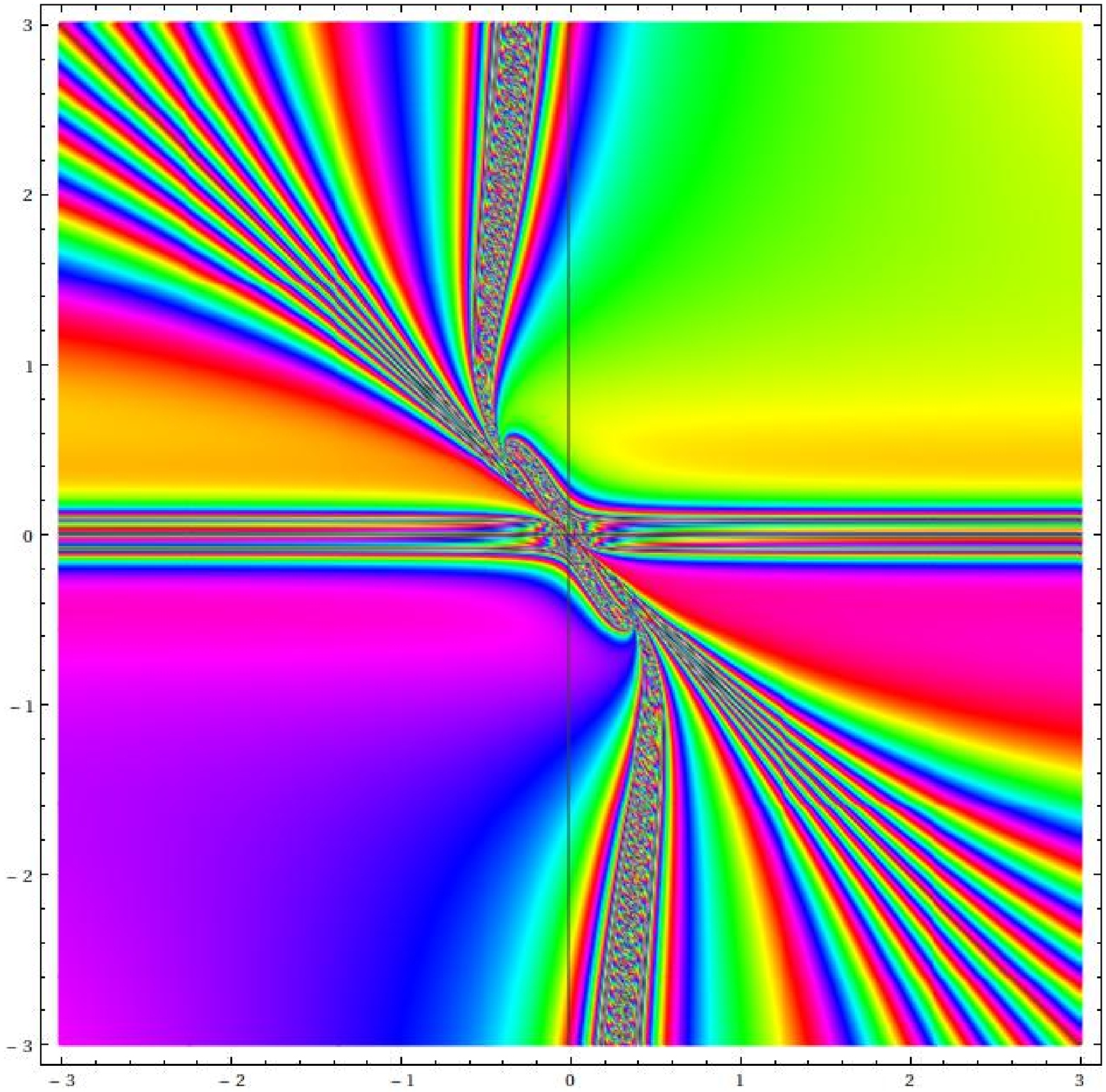} &
\includegraphics[width=0.45\textwidth]{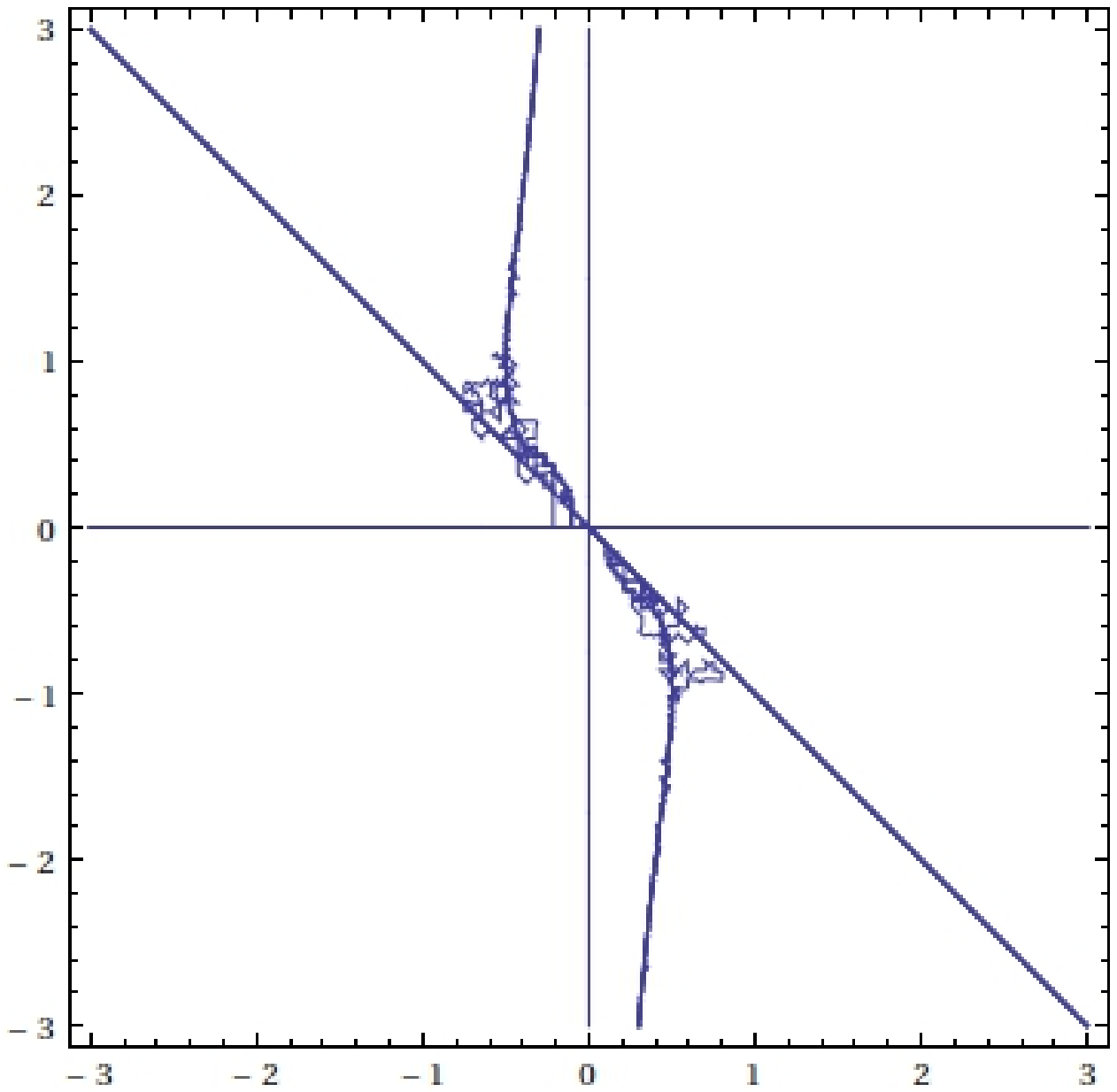}\\
\includegraphics[width=0.45\textwidth]{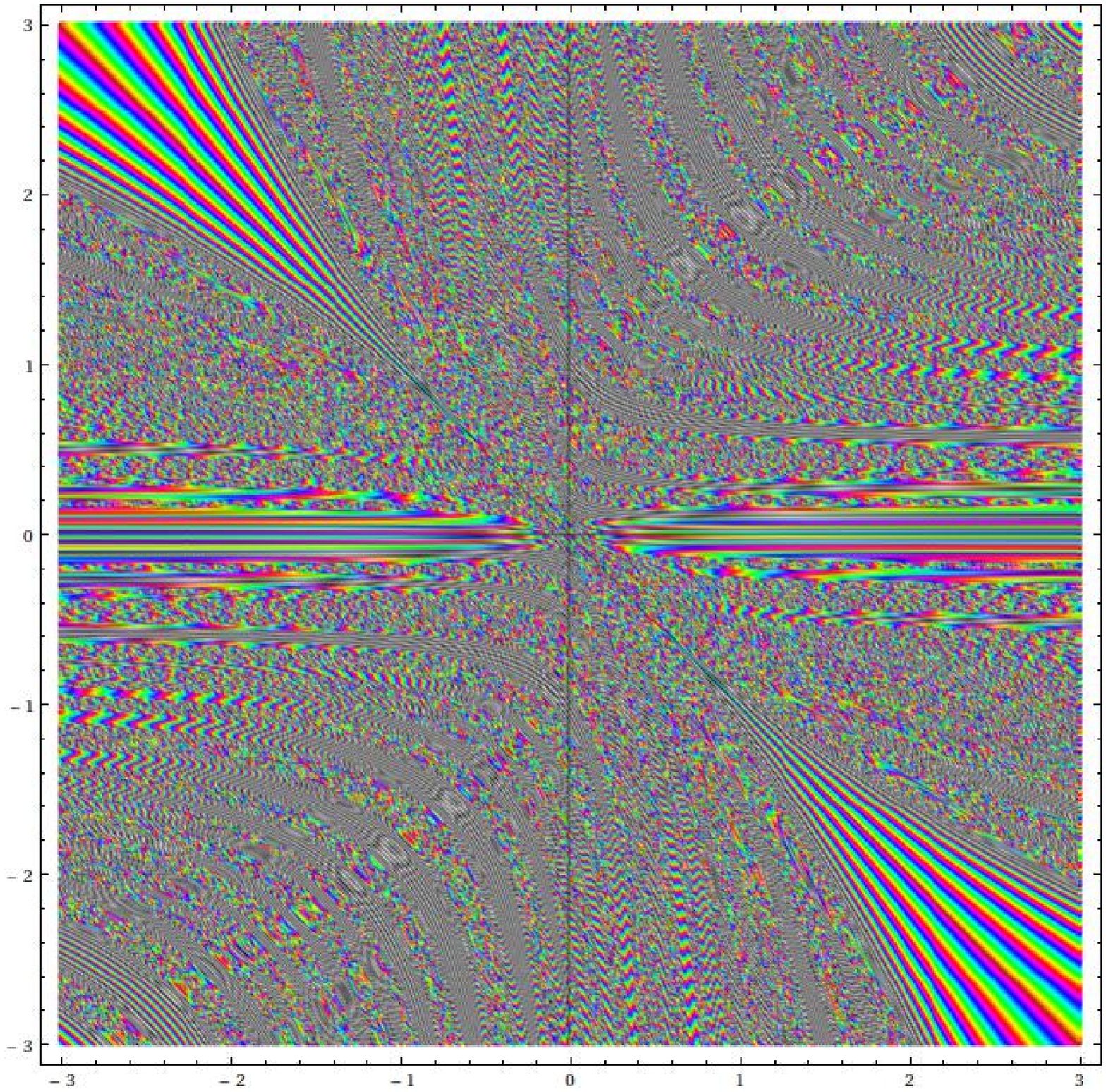} &
\includegraphics[width=0.45\textwidth]{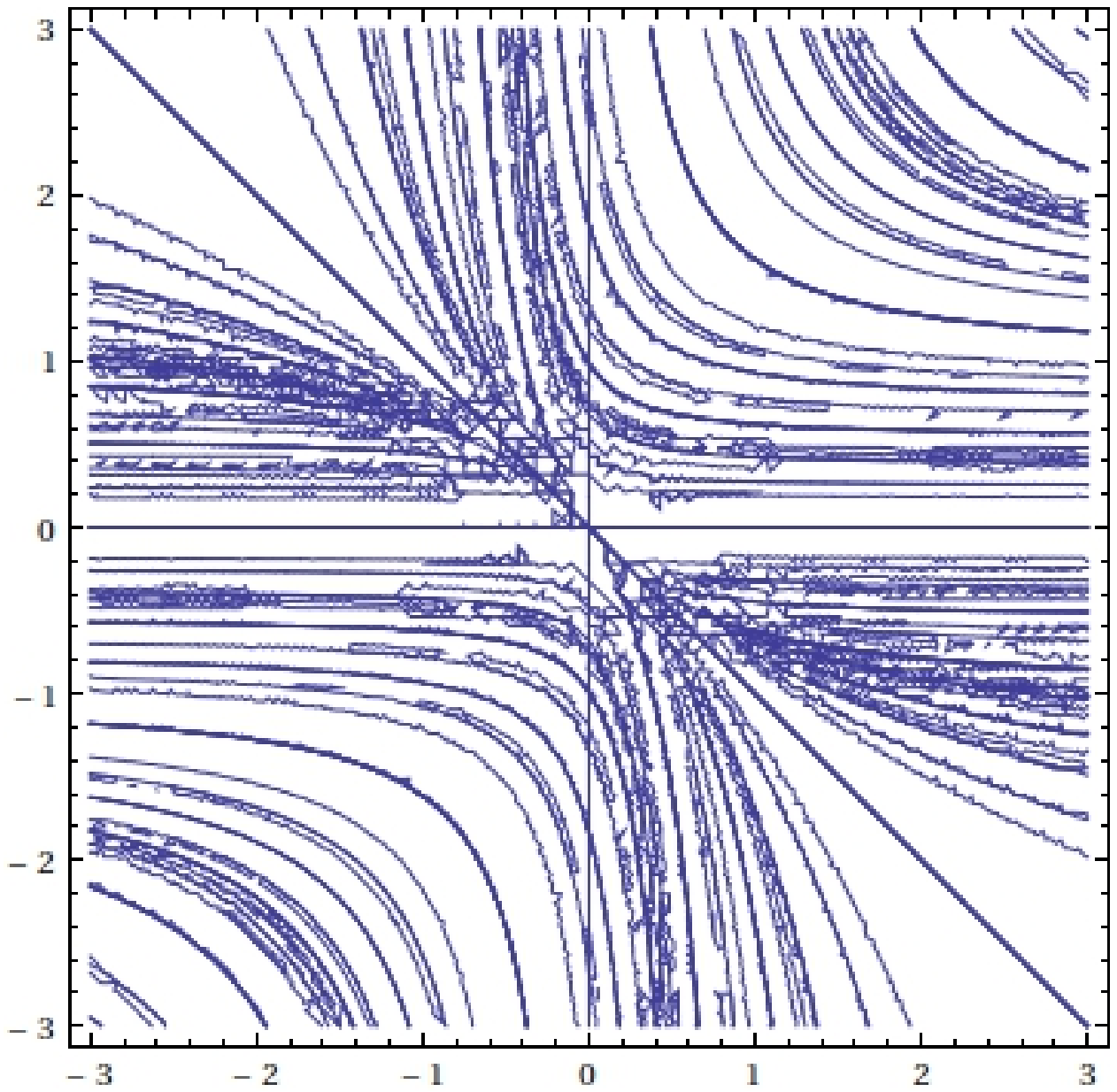}
\end{array}$
\end{center}
\caption{Two representations of the FS of DE (\ref{eq:DeVaultGalminasJanowskiLadas}). On the right the first $10$ forbidden curves are drawn. On the left we color after the value of the tenth iterate of each point, reserving gray for points of $\mathcal F$. From up to down, the values of the parameters $(A,B)$ are the following: $(1,5)$, $(1,1)$, $(-1,-1)$. Numerical experiments show also that when $(A,B)$ equals to $(1,-1)$ and $(-1,1)$ respectively, vertical symmetries of the second and third row are obtained.}
\label{fg:DeVaultGalminasJanowskiLadas}
\end{figure}

In a similar way, we can estimate the FS of 
\begin{equation}
x_{n+1}=\frac A {x_n}+ \frac B{x_{n-1}}
 \label{eq:DeVaultGalminasJanowskiLadas}
\end{equation}
for differents values of parameters $A$ and $B$. See figure \ref{fg:DeVaultGalminasJanowskiLadas} and open problem \ref{op:DeVaultGalminasJanowskiLadas}.
\section{A list of open problems}
\begin{OP}
Clarify which kind of topological objects can be obtained in the second order Riccati DE and generalize to the Riccati DE of order $k$ (\ref{eq:RiccatiOrderk}).
\label{op:QualitativeDescriptionRiccati} 
\end{OP}
\begin{conj}
DE (\ref{eq:RiccatiOrderk}) is globally periodic if and only if $\mathcal F$ is a finite collection of hypersurfaces in $\R^k$.
 \label{conj:RiccatiGloballyPeriodic}
\end{conj}
\begin{OP}
Describe the FS of RDE (\ref{eq:AboZeid_generalized}).
  \label{op:AboZeid_generalized}
\end{OP}
\begin{OP}
Describe the FS of RDEs (\ref{eq:BajoLiz2011_generalized1}) and (\ref{eq:BajoLiz2011_generalized2}).
 \label{op:BajoLiz2011_generalized}
\end{OP}
\begin{OP}
Describe the FS of RDEs (\ref{eq:Mcgrath2006_generalized}).
 \label{op:Mcgrath2006_generalized}
\end{OP}
\begin{OP}
Let $T(x)=\frac{\alpha x+\beta}{\gamma x + \delta}$. Describe the FS of the RDEs obtained by applying the change $x_n=T(y_n)$ to DE (\ref{op:Rhouma_multiplicative}).
\label{op:Rhouma_multiplicative}
\end{OP}
\begin{OP}
Let $T(x)=\frac{\alpha x+\beta}{\gamma x + \delta}$. Describe the FS of the RDEs obtained by applying the change $x_n=T(y_n)$ to DE (\ref{op:Rhouma_linear}).
\label{op:Rhouma_linear}
\end{OP}
\begin{OP}
Let $x_{n+1}=G(x_n,\ldots,x_{n-k})$ a DE of order $k+1$ such that its closed form is known. Let $x_n=H(y_{n+1},\ldots,y_{n-k})$ a change of variables that can be rewriten as a nonautonomous DE $y_{n+1}=\tilde{H}(y_n,\ldots,y_{n-k})$ with coefficients depending of $x_n$ and such that its closed form solution is also known. Determine the forbidden set of the DE obtained by applying the change of variables to $x_{n+1}=G(x_n,\ldots,x_{n-k})$.
\label{op:Rhouma_generalized}
\end{OP}
\begin{OP}
 Describe the FSs of DEs (\ref{eq:Palladino_generalized_equation}) and (\ref{eq:Palladino_generalized_equation_order_k}) having algebraic invariants.
 \label{op:Palladino_generalized}
\end{OP}
\begin{OP}
Let $k$ and $l$ be distint natural numbers greater than $0$. Determine the closed form solution and the forbidden set of DEs admiting an invariant of the form (\ref{eq:AghajaniShouli_invariant}).
 \label{op:AghajaniShouli}
\end{OP}

\begin{OP}
For every RDE of order $2$, determine a region $A\subset\R^2$ such that the forbidden curves in $A$ can be given in explicit form, and determine them.
 \label{op:ExplicitCurves}
\end{OP}
\begin{OP}
Find the symbolic description of $\mathcal F_{\C}$ and $\mathcal F_{\R}$ for the families of RDEs (\ref{eq:InverseParabola}) and (\ref{eq:InverseLogistics}).
 \label{op:SymbolicDescription}
\end{OP}
\begin{OP}[\cite{RubioMassegu2009}, open problem 1]
 Let $F:D\subseteq \R^k\rightarrow \R^k$ an almost local diffeomorphism such that $F^{-1}(\R^k\setminus D)$ has Lebesgue measure zero. Suppose that equation (\ref{eq:RubioMassegu_Map}) is globally periodic. Obtain sufficient conditions in order that $\mathcal F$ is closed.
  \label{op:RubioMassegu}
 \end{OP}
\begin{conj}[\cite{RubioMassegu2009}, conjecture 2]
 Let $f:D\rightarrow\C$ be a rational function and $D$ its natural domain. Assume that $f$ depends effectively on its last variable. Then the good set of (\ref{eq:RubioMassegu_DE}) is open if, and only if, the equation is globally periodic.
\label{conj:RubioMassegu}
 \end{conj}
\begin{OP}[\cite{RubioMassegu2007OnTheExistence}, open problem 1.]
To obtain necessary conditions in order to have that an equilibrium point of an autonomous difference equation belongs to the interior of the good set.
\label{op:RubioMassegu_EquilibriumInTheInterior}
\end{OP}
\begin{OP}[\cite{RubioMassegu2007OnTheExistence}, open problem 2.]
To obtain necessary conditions in order to have that the good set of a difference equation is open.
\label{op:RubioMassegu_OpenGoodSet}
\end{OP}
\begin{OP}
Let $f:\R\rightarrow \R$ be a monotonic map with pole. Complete the results of theorem \ref{th:Sedaghat_MonotonicMapWithPole}, in the sense of describing qualitatively the FS of DE $x_{n+1}=f(x_n)$ when $f$ is an increasing map with no fixed point. In particular, answer to the following questions:
\begin{enumerate}
 \item In which cases $\mathcal F$ is finite?
 \item Can be $\mathcal F$ dense outside a compact set? (See \cite{Sedaghat2000}).
\end{enumerate}
\label{op:Sedaghat_DenseFS}
\end{OP}
\begin{OP}[\cite{2000_DeVaultGalminasJanowskiLadas}]
 Determine the FS of DE (\ref{eq:DeVaultGalminasJanowskiLadas}).
 \label{op:DeVaultGalminasJanowskiLadas}
\end{OP}
\begin{OP}
Determine the FS of DE (\ref{eq:Sedaghat_Qualitative_1_generalized}) by applying theorem \ref{th:Sedaghat_MonotonicMapWithPole} in terms of $a$, $b$, $c$ and $p$.
 \label{op:Sedaghat_Qualitative_1_generalized}
\end{OP}

As we have seen throughout the paper, there are a lot of DEs with closed form solution and explicit forbidden set. There are many ways to generalize them and therefore the literature on this subject is growing with more and more works devoted to that kind of DE. Although those works have the value to enrich the set of known examples, may be it is time to put some order in the field. We propose to elaborate a database of DEs and SDEs including, among others, the following topics: the forbidden and good sets, the asymptotic behavior, the boundedness of solutions, the relationships with other DE or SDE and the applied models related with them.
%
\section{Acknowledgements}
This paper has been partially supported by Grant MTM2008-03679 from Ministerio de Ciencia e Innovac\'on (Spain), Project 08667/PI-08 Fundaci\'on S\'eneca de la Comunidad Aut\'onoma de la Regi\'on de Murcia (Spain).\\
We thank the anonymous referees for their suggestions and comments.

\bibliographystyle{plain} 

\end{document}